\documentclass[10pt]{article}

\usepackage{amssymb,amsmath,latexsym}
\voffset
-1.5cm

\usepackage{a4,amscd,graphicx,epsfig}
\everymath{\displaystyle}

\newtheorem{teo}{Theorem}[section]
\newtheorem{lem}[teo]{Lemma}
\newtheorem{cor}[teo]{Corollary}

\newtheorem{prop}[teo]{Proposition}
\newtheorem{defi}[teo]{Definition}

\newtheorem{remark}[teo]{Remark}

\newcommand{\mr}{\mathbb{R}}
\newcommand{\mc}{\mathbb{C}}
\newcommand{\mz}{\mathbb{Z}}
\newcommand{\mh}{\mathbb{H}}

\newcommand{\Bb}{{\mathcal B}}
\newcommand{\Cc}{{\mathcal C}}
\newcommand{\Dd}{{\mathcal D}}

\newcommand{\Ff}{{\mathcal F}}
\newcommand{\Gg}{{\mathcal G}}
\newcommand{\Hh}{{\mathcal H}}
\newcommand{\Ii}{{\mathcal I}}

\newcommand{\Kk}{{\mathcal K}}
\newcommand{\Ll}{{\mathcal L}}
\newcommand{\Mm}{{\mathcal M}}

\newcommand{\Rr}{{\mathcal R}}

\newcommand{\Tt}{{\mathcal T}}

\newcommand{\T}{{\mathbb T}}
\newcommand{\C}{{\mathbb C}}

\newcommand{\Z}{{\mathbb Z}}
\newcommand{\R}{{\mathbb R}}
\newcommand{\N}{{\mathbb N}}
\newcommand{\A}{{\mathbb A}}

\title{3D Quantum Hyperbolic Field Theory}

\author {St\'ephane Baseilhac$^1$ and Riccardo Benedetti$^2$}

\date {}

\begin{document}

\maketitle

\vspace{0.7cm}

\noindent
$^1$ Universit\'e de Grenoble I, Institut Joseph Fourier, UMR CNRS
5582, 100 rue des Maths, B.P. 74, F-38402 St Martin d'Heres Cedex.
Email: baseilha@ujf-grenoble.fr

\smallskip

\noindent 
$^2$ Dipartimento di Matematica, Universit\`a di Pisa, Via F.
Buonarroti, 2, I-56127 PISA. Email: benedett@dm.unipi.it.

\vspace{1cm}
\begin{abstract}
\noindent We construct a new family of exact quantum field
theories modeled on hyperbolic geometry, called {\it quantum
  hyperbolic field theories} (QHFTs). All the QHFTs are defined for a same
$(2+1)$-bordism category, based on the set of compact oriented
$3$-manifolds $Y$, equipped with properly embedded framed links
$L_\Ff$ and with flat connections $\rho$ of principal
$PSL(2,\C)$-bundles over $Y \setminus L_\Ff$, with arbitrary holonomy
at the link meridians. The QHFTs generalize our previous works
\cite{BB0,BB1,BB2} on volumes, Chern-Simons invariants and quantum
hyperbolic invariants of $PSL(2,\C)$-characters (i.e. conjugacy
classes of $PSL(2,\C)$-valued representations of the fundamental
group) of closed $3$-manifolds. A main part of the paper consists in
specifying the {\it marked surfaces} that make the objects of the
bordism category. This marking includes the introduction of adequate
parameters for the space of all $PSL(2,\C)$-characters of a punctured
surface, which is the fundamental {\it phase space} of the QHFTs. Each
QHFT associates to a triple $(Y,L_\Ff,\rho)$ as above with marked
boundary components a tensor, which is generically holomorphic w.r.t.
the parameters for the restriction of $\rho$ to the punctured boundary
$\partial Y \setminus
L_\Ff$. As a first application, we get new numerical invariants of
$3$-manifolds, such as Chern-Simons invariants of
$PSL(2,\mc)$-characters of arbitrary link complements, or quantum
invariants of compact hyperbolic cone manifolds. Another application
is the construction, for any $PSL(2,\mc)$-character of a surface of
finite topological type, of new conjugacy classes of linear
representations of the mapping class group. Finally, we discuss some
evidences showing that the QHFTs are pertinent to 3D gravity.
\end{abstract}


\section{Introduction}\label{DQFTINTRO}

This paper is the continuation of our previous works
\cite{BB0,BB1,BB2}. It completes the construction of a new family of
3D-quantum field theories (QFTs), whose fundamental ingredients are,
on the combinatorial side, suitably structured families of hyperbolic
ideal tetrahedra, and, on the functional side, matrix versions of the
dilogarithm functions.  As both play a central role in
$3$-dimensional hyperbolic geometry and volume computations, we call
these QFTs {\it quantum hyperbolic field
  theories} (QHFTs).

More precisely, having as model Atiyah's formalization of the
{\it topological quantum field theories} (TQFTs) \cite{At,T}, we use
the terms ``3D-quantum field theory'' as synonimous of: \smallskip

\noindent {\it Monoidal functor from some $(2+1)$-bordism category to the
tensorial category of complex linear spaces.}
\smallskip

\noindent Recall that it means a correspondence (a ``representation'')
sending {\it marked} surfaces to complex linear spaces, and
3-dimensional manifolds to linear morphisms between the linear spaces
associated to the marked boundary components. This correspondence maps
the gluing of bordisms to the morphism composition, and respects as
well certain tensor products in both categories. How the marking of
surfaces is rich does reflect in how the part of the QFT supported by
the product bordisms is non trivial; at least it should include
projective representations of the appropriate mapping class groups.

\noindent In the case of TQFTs one uses essentially the bare
topological bordism category, but this set up can be extended to
bordism categories supported by suitably ``equipped'' $3$-manifolds.
The family of QHFTs that we construct in this paper is defined for a
$(2+1)$-bordism category based on oriented compact $3$-manifolds $Y$,
which are equipped with properly embedded framed links $L_\Ff$ and
with flat connections (up to gauge equivalence) on principal
$PSL(2,\C)$-bundles over $Y \setminus L_\Ff$ (i.e. with
$PSL(2,\C)$-characters of $Y \setminus L_\Ff$), having {\it arbitrary}
holonomy at the meridians of the link components. Recall that
$PSL(2,\mc)$ can be identified with the group of direct isometries of
the hyperbolic $3$-space.

\smallskip

In fact, the main themes of Thurston's geometrization
program, such as the Teichmuller spaces with the action of the modular
groups, the hyperbolic volume, or more generally numerical invariants of $PSL(2,\C)$-characters of
$3$-manifolds like the volume or the Chern-Simons invariants, have
quantum analogs which are contained in the QHFTs, in the sense that
they are completely described by these field theories. So the
differential geometry of these central classical objects should hopefully
reflect deeply in the QHFTs.

\noindent For instance, we have defined in \cite{BB1,BB2} new families of
complex valued (up to a determined phase ambiguity) quantum invariants
for cusped hyperbolic $3$-manifolds, and for arbitrary triples
$(W,L,\rho)$, where $W$ is a compact closed oriented $3$-manifold, $L$
is an unframed link in $W$, and $\rho$ is a $PSL(2,\C)$-character
defined on the whole of $W$ (hence it is trivial at the link
meridians). We also obtained new simplicial formulas for the volume
and the Chern-Simons invariants of such characters. As explained in
Section \ref{TENSOR} of the present paper, all these invariants are
specific QHFT {\it partition functions} (using the classical
terminology of the physics litterature), or variations of them. The
celebrated Volume Conjectures state precise relationships between
the 'semi-classical limit' of the quantum invariants and the volume
and the Chern-Simons invariants, when the manifold admits a hyperbolic
structure and the $PSL(2,\C)$-character is the hyperbolic holonomy. We
refer to Section 5 of \cite{BB1}, Sections 6-7 of \cite{BB2}, and to
\cite{B1} for details on these conjectures.

\smallskip

In this paper, we are mainly concerned with the extension of the heavy
apparatus of combinatorial structures underlying the simplicial
formulas of these invariants to
manifolds {\it with marked boundary components}. In particular, we
introduce several parameter spaces for the basic {\it phase space} of
the theory, which is the space $\Rr(g,r)$ of $PSL(2,\C)$-characters of
punctured surfaces $\Sigma_{g,r}$. Each parameter space is an
algebraic variety that fibers over $\Rr(g,r)$ and admits a
filtration, with, for instance, different stratas being bundles over
the Teichmuller spaces of hyperbolic closed surfaces, punctured
surfaces, or surfaces with totally geodesic boundary.

\noindent The parameter's construction is self-contained, has a purely
$3$-dimensional interpretation, and is naturally adapted to the QHFTs.
However, it is remarkable that the $\Ii\partial$-{\it parameter}
spaces (on which the QHFT morphisms are eventually defined) are highly
reminiscent of the well-known shear-bend coordinates for pleated
hyperbolic surfaces with punctures. On another hand, the full QHFT
marking of surfaces, including these spaces of parameters, encode the
irreducible representations of a ``quantum moduli space'' of
$PSL(2,\mc)$-characters of punctured surfaces, very similar to the
Kashaev or (exponential version of the) Chekhov-Fock quantum
Teichmuller spaces (see \cite{K3}, \cite{C-F}). We plan to investigate
both facts in a separate paper.

\smallskip

The QHFTs are \emph{exact} (in principle, every QHFT morphism can be
explicitely computed), {\it finite dimensional} (i.e. the linear space
associated to any marked surface is finite dimensional), and {\it
  hermitian}. They form a family indexed by the odd positive integers
$N \geq 1$. Each QHFT associates to a triple $(Y,L_\Ff,\rho)$ as above
with marked boundary a tensor, which is holomorphic (up to a
determined ambiguity) on a dense subset of the $\Ii\partial$-parameter
space for the $PSL(2,\mc)$-characters of the (punctured) boundary of
$Y \setminus L_\Ff$. In particular, we get tensor valued holomorphic
functions on the bundles of $\Ii\partial$-parameters over
the Teichmuller spaces for the boundary components.

\noindent For instance, in
the case of product bordisms, by varying the marking simultaneously on
both boundary components, letting fixed the character, these tensors
define conjugacy classes of projective reprensentations of the mapping
class groups\footnote{The resulting invariants of surface
  diffeomorphisms should be closely related to those obtained recently by
  Bonahon-Liu \cite{Bon-Liu}.}. Also, in the ``classical'' case when
$N=1$, these tensors are just scalars.  By extending the results of
\cite{BB2}, Section 6 (which hold for closed manifolds or cusped
hyperbolic manifolds), they can be interpreted as the evaluation of a
second Cheeger-Chern-Simons class for manifolds $Y$ with {\it
  parametrized boundary}, i.e. as ${\rm CS}(\rho) + \sqrt{-1}{\rm
  Vol}(\rho)$, where $\rho$ is a $PSL(2,\mc)$-character of $Y$, and
Vol and CS are respectively a volume and Chern-Simons invariant of
$\rho$ on the {\it marked} bordism $Y$. Finally, the QHFT morphisms
for triples $(W,L_\Ff,\rho)$, where $W$ is a closed manifold, are
always scalars. It is an open problem to understand their dependence
w.r.t. the framing of $L_\Ff$, as well as the relationship between
them and the invariants of cusped hyperbolic $3$-manifolds defined in
\cite{BB2} (compare with Section \ref{VARIATION}).

\smallskip

The ultimate building blocks of the QHFTs are so called {\it matrix
  dilogarithms}, which are determined automorphisms $\Rr_N$, $N \geq
1$, of $\C^N\otimes \C^N$ associated to hyperbolic ideal tetrahedra
equipped with an elaborated extra-decoration, and that satisfy certain
fundamental {\it five term identities}. The matrix dilogarithms have
been introduced, formalized and widely studied in \cite{BB2}. They are
``quantum'' versions of the classical dilogarithm functions (see
Section 8 of \cite{BB2} and the references therein, and \cite{B1}).
 
\smallskip

The above deep interaction between classical objects coming from
differential geometry and analysis, and quantum algebraic objects, is
not the only motivation for studying the QHFTs. Another one is the
fact that the whole family of QHFTs forms a unified theory that could
be understood as a finite regularization of quantum 3D gravity. This
is discussed in Subsection \ref{3Dgrav} below. Before that, we describe the
content of the paper in the next Subsection \ref{descpaper}.

\smallskip

Let us conclude this introduction by noting that Turaev has formalized
in \cite{Tu2} a notion of Homotopic QFT (HQFT), which provides a
general framework for QFTs based on $3$-dimensional cobordisms
equipped with a representation of their fundamental group in a fixed
group G. It is tempting to look at the QHFTs as examples of HQFTs for
G=PSL(2,C), but some key points in our construction show that this
cannot be exactly the case. For instance, there is the needed ``link
fixing'', which implies that we have to use punctured surfaces as
objects. However the relationship between the QHFTs and Turaev's
HQFTs certainly deserves further investigations.

\subsection{Description of the paper} \label{descpaper}

The adequate sets of parameters for the spaces of
$PSL(2,\C)$-characters of punctured surfaces are developped in Section
\ref{PARAM}. This section is self-contained; only some notions, which
we recall when needed, are taken from \cite{BB1,BB2}.

After some preliminaries on the variety of $PSL(2,\C)$-characters, in
Subsection \ref{e-triang}, we introduce for any closed compact surface
$S$ with a finite set $V=\{p_i\}_{i=1}^r$ of framed (i.e. with a fixed
segment $l_i$ in a disk neighborhood) marked points a notion of {\it
  efficient triangulation} for the surface with boundary $\textstyle F
= S \setminus \coprod_i D_i$, where $D_i$ is a small open disk with
center on $l_i$. These triangulations have two main features. First,
they are naturally adapted to the $3$-dimensional machinery of
$\Ii$-cusps developped in Section \ref{BORDCAT}. Second, they allow to
define charts for the {\it whole} space of $PSL(2,\C)$-characters of
$F$, for any kind of boundary holonomies. Namely, given any efficient
triangulation $T$ of $F$ with a suitable system $b$ of orientation of
the edges called {\it branching}, we produce bundles of cocycle
$\Dd$-{\it parameters} (Subsection \ref{Dpar})
$$\ Z(T,b) \longrightarrow R(g,r)$$
and $\Ii\partial$-{\it parameters} (Subsection \ref{Ipar})
$$W(T,b) \longrightarrow R_I(g,r)$$
where $R(g,r)$ denotes the space
of $PSL(2,\C)$-characters of $F$ and $R(g,r)_I$ is a Zariski open
subset of $R(g,r)$. The total spaces of these bundles are algebraic
varieties, which admit partitions into subbundles, according to the
type of the holonomies at the boundary components of $F$ (trivial,
parabolic, or else). The total spaces of these subbundles form a
filtration by the dimension. Examples are the parameters for the
Teichmuller space of the closed surface with marked points $(S,V)$ (in
the deepest part of the filtration), and the parameters for the
Teichmuller space of finite area complete hyperbolic metrics on the
punctured surface $S \setminus V$ (in the quasi-regular part of the
filtration).

\noindent The QHFT morphisms eventually depend on the $\Ii\partial$-parameters,
but the cocycle $\Dd$-parameters come at first naturally from the
combinatorial presentation of manifolds we need for all the
construction. Namely, they are the matrix entries of
$PSL(2,\mc)$-valued $1$-cocycles on efficient triangulations that
represent the elements of $R(g,r)$, with a preferred ``normalized''
form into which {\it any} cocycle can be put by conjugation. The
$\Ii\partial$-parameters are products of cross-ratio moduli of certain
families of decorated hyperbolic ideal tetrahedra. These hyperbolic
tetrahedra are associated to the $3$-simplices of standard
triangulated cylinders over $F$, by suitably extending the above
normalized cocycles to the cylinders, and then using an {\it
  idealization procedure} reminiscent of the construction of
piecewise-straight developping maps for geometric structures on
$3$-manifolds. The $\Ii\partial$-parameters can also be viewed as
describing representations of the groupoid of paths {\it transverse}
to the given efficient triangulation (see the end of Subsection
\ref{Ipar}).

\noindent For the convenience of the reader, we present in an Appendix the
relationship between the cocycle $\Dd$-parameters and the
Kashaev-Penner coordinates for the moduli space of irreducible
$PSL(2,\mc)$-characters on punctured surfaces, with parabolic
holonomies at the punctures (this includes the Teichmuller space of $S
\setminus V$).

\smallskip

In Section \ref{BORDCAT} we define the QHFT $(2+1)$-bordism category,
based on triples $(Y,L_\Ff,\rho)$ (Subsection \ref{bord+}). In
  particular we carefully describe the {\it marked surfaces} that make
  the (elementary) objects of the category. In fact, it is more
  convenient to deal with an equivalent category based on
  $3$-manifolds with corners, at the intersection of a closed tubular
  neighborhood of $L_\Ff$ with $\partial Y$. This makes the inclusion
  of the results of Section \ref{PARAM} immediate, as the marking of
  surfaces shall incorporate the phase space parameters.  Finally, we
  introduce in Subsection \ref{CUSP} the notion of $\Ii$-{\it cusps}:
  these are standard forms to represent pairs $(U(L_\Ff),\rho_\vert)$,
  where $U(L_\Ff)$ is a closed tubular neighborhood of $L_\Ff$ and
  $\rho_\vert$ is the restriction of $\rho$ onto it, as the gluing of
  suitably decorated hyperbolic ideal tetrahedra. This notion is
  essential to the QHFTs, for instance to obtain numerical invariants
  for closed manifolds $Y$, in the case when the character $\rho$ is
  non trivial at the meridians of $L_\Ff$.

\smallskip

In Section \ref{TENSOR} we quickly review the matrix dilogarithms (the
explicit formulas are given in an Appendix) and we define the
tensors that represent the bordisms. This completes the construction of
the QHFTs. We heavily refer to the notions and results developped
in \cite {BB0,BB1, BB2}, avoiding too many unecessary repetitions, and
pointing out the substantial new achievements. 

In particular, we discuss the corresponding QHFT partition functions
for {\it closed} manifolds $Y$, and a variation of the QHFT
construction (Subsection \ref{PARTF}). Namely, we consider a more restricted bordism category such that
the associated partition functions for triples $(W,L,\rho)$,
where $W$ is closed, $L$ is an {\it unframed} link, and $\rho$ is {\it
  trivial} at the link meridians, coincide with the dilogarithmic
invariants already defined in \cite{BB1, BB2}. We set the
relationship between the two kinds of partition functions that are
available for triples $(W,L_\Ff,\rho)$, for such a special $\rho$.

We consider also the part of the QHFTs supported by the trivial
(product) bordisms (Subsection \ref{MODREP}). As mentioned in the
Introduction, it contains interesting conjugacy classes of linear
representations of the mapping class groups of punctured surface
(defined up to a determined phase ambiguity).

Finally we indicate in Subsection \ref{VARIATION} a so called {\it universal
  QHFT environment}, that is the most general set up where our
constructions {\it formally} makes sense. The specialized QHFTs
previously constructed naturally map into this universal environment,
and {\it have a clear intrinsic topological/geometric meaning}. This
suggests the possibility of other meaningful specializations of the
universal QHFT environment.

\subsection{QHFT and 3D gravity}\label{3Dgrav}
There are some evidences that the QHFTs are pertinent to 3D gravity
(see \cite{Be} for further comments on this points). Thanks to the
$3$-dimensional peculiar fact that the Ricci curvature tensor
completely determines the Riemann curvature tensor, classical 3D
(pure) gravity concerns the study of Riemannian or Lorentzian
$3$-manifolds of {\it constant curvature}. The sign of the curvature
coincides with the sign of the cosmological constant for the theory.
We stipulate that all manifolds are {\it oriented} and that the
Lorentzian space-times are also {\it time-oriented}. We also include
in the picture the presence of {\it world lines} of ``particles''; the
singularies of the metric are concentrated along these lines. Typical
examples are the {\it cone manifolds} of constant curvature with a
properly embedded link as cone locus, where the cone angles reflect
the ``mass'' of the particles. In the Lorentzian case we also require
that the world lines are of causal type (see e.g. \cite {BG2}).

\noindent 
The hyperbolic $3$-manifolds are the classical solutions of Riemannian
(sometimes called {\it Euclidean}) 3D gravity with (normalized) negative
cosmological constant. {\it Geometrically finite}
hyperbolic $3$-manifold, or more generally topologically {\it tame}
ones, possibly with links of concentrated singularities (hence with
non necessarily trivial holonomy at the link meridians), give
fundamental examples of supports for the QHFT bordism category. Here it is
understood that these manifolds are equipped with the holonomies of
the hyperbolic structures (remind that $PSL(2,\C)$ is identified with
${\rm Isom}^+(\mh^3)$, the group of direct isometries of the
hyperbolic $3$-space). Moreover, for compact hyperbolic $3$-manifolds
$W$, a deep {\it volume rigidity} result (see
e.g. \cite{Dun, F2}) tells
us that a volume function is well defined on the space $\Rr(W)$ of
conjugacy classes of $PSL(2,\C)$-valued representations of $\pi_1(W)$,
and that, if $\rho$ is the holonomy of a hyperbolic structure $h$ on
$W$, then: \smallskip

(1) $ {\rm Vol}(\rho) = {\rm Vol}(W,h)$;
\smallskip

(2) $\rho$ is the unique maximum of the volume function.
\medskip

\noindent (With some technical complication this result holds also for {\it
  cusped} manifolds, i.e. for non compact finite volume complete
hyperbolic $3$-manifolds). This geometric result is strictly related
to Euclidean 3D gravity with negative cosmological constant, when
formulated in terms of a Chern-Simons type action for the so called
``new variables'', which are the connections on principal
$PSL(2,\C)$-bundles, instead of the metrics and the framings (see
\cite{W}). The ``constraint'' equations for this action imply that the
phase space of the theory becomes the space of {\it flat}
$sl(2,\mc)$-connections (up to gauge equivalence). In fact, the
Chern-Simons action for flat $sl(2,\mc)$-connections equals a constant
times $ {\rm CS}(\rho)+ i{\rm Vol}(\rho)$, where ${\rm CS}(\rho)$
denotes the Chern-Simons invariant of the flat connection $\rho$. This
is a natural complexification of the above volume function on
$\Rr(W)$, and hyperbolic manifolds, the classical solutions of
Riemannian 3D gravity, maximize the norm of $\exp((1/2i\pi)({\rm
  CS}(\rho)+ i{\rm Vol}(\rho)))$. It is a fact (see \cite{BB2} and
Section 5 of this paper) that the ``classical member'' QHFT$_1$ of the
family actually computes this exponentiated classical complex action
for pairs $(W,\rho)$. Moreover, different so called ``Volume
Conjectures'' should identify QHFT$_1$ with QHFT$_\infty$, i.e. the
``classical limit'' of the ``quantum'' theories QHFT$_N$, $N>1$, when
$N\to \infty$ (see \cite {BB1, BB2} or \cite{B1} for a discussion on
this point). Remind also that, in the particular case of a link $L$ in
$S^3$ equipped with the {\it trivial} flat bundle, QHFT$_N$, $N>1$,
computes the Kashaev's \cite{K1} invariant $<L>_N$, later identified
by Murakami-Murakami \cite{MM} with $J_N(L)(\exp(2\pi i/N)$, where
$J_N$ denotes a suitably normalized colored Jones invariant.

Another intriguing fact is that bordisms supported by hyperbolic
$3$-manifolds are {\it not only} pertinent to Euclidean 3D gravity
with negative cosmological constant. This claim comes from the
following few facts; we refer to \cite{M, BG1, Bo, BBo} for the
details and more articulated statements. We recall, for example, that
geometrically finite hyperbolic $3$-manifolds with incompressible ends
of infinite volume can be concretely interpreted as {\it interactions}
between {\it Lorentzian} space-times of arbitrarily fixed constant
curvature.  More precisely, we can canonically associate to every end
of such a hyperbolic $3$-manifold a domain of dependence of a compact
Cauchy surface, of {\it arbitrarily fixed} constant curvature
$\kappa$. A key point is that these Lorentzian space-times,
independently on their constant curvature, share the same ``parameter
space'' $T_g\times \Mm\Ll_g$, where $T_g$ denotes the Teichm\"uller
space of hyperbolic structure on a fixed surface $S$ of genus $g\geq
2$, and $\Mm\Ll_g$ is the space of {\it measured geodesic laminations}
(see e.g. \cite{Ep}) on these hyperbolic surfaces. Moreover this
is also the parameter space of {\it projective structures} on $S$
\cite{Th}. The holonomy of the projective structure related to a
hyperbolic end as above is just the restriction of the holonomy of the
hyperbolic $3$-manifold. Moreover these space-times have a very
explicit geometric description. In particular, they admit a canonical
{\it cosmological time}: the proper time that every event has been in
existence and that coincides with its {\it finite} Lorentz distance
from the {\it initial singularity}. This initial singularity has a
rich geometry (``dual'' to the geodesic lamination) which is the
``past limit'' in an appropriate sense of the geometry of the level
surfaces of the cosmological time. Every such a level surface is a
Cauchy surface.  Moreover, when $\kappa \leq 0$, there is a canonical
{\it Wick rotation} directed by the gradient of the cosmological time
(which is in general a $C^0$ vector field) that converts the future of
a determined level surface into the whole associated hyperbolic end.
Remind that the very basic example of Wick rotation directed by the
field $\partial/\partial x_3$ converts the Minkowski metric on $\R^3$
with signature $(+,+,-)$ into the Euclidean metric; sometimes one
refers to it as ``passing to the imaginary time''. Wick rotation is a
basic procedure for interplaying Riemannian and Lorentzian geometry,
including the {\it global causality} of Lorentzian space-times.

In a sense, this behaviour confirms the intuition at page 72 of
\cite{W}: these Lorentzian space-times should be considered not really
as ``space-times'', but rather as mere ``world sheets''; hence it does
not really make sense to ask about their curvature. The latter is
matter of a ``universe'' where they should be embedded.  The above
considerations show in particular that hyperbolic universes can
concretely realize the {\it changes of topology} of these world
sheets in a purely {\it classical} 3D gravity set up,
providing that we avoid any (somewhat missleading) separation in
different sectors, accordingly to the metric signature and the sign of
the cosmological constant.


\section{ Phase space parameters}\label{PARAM}

To orient the boundary $\partial Y$ of any oriented $n$-manifold $Y$,
we adopt the convention: {\it last is the ingoing normal}. 

\smallskip
 
For every $(g,r)\in \N \times \N$, such that $g\geq 0$, $r>0$, and
$r>2$ if $g=0$, we fix a compact closed oriented {\it base} surface
$S=S_g$ of genus $g$, with a set $V=V_{g,r}=
\{v_1,\dots \ , v_r\}$ of $r$ marked points. Our basic ``phase space'' is
$$\Rr(g,r)= {\rm Hom}(\pi, PSL(2,\C))/PSL(2,\C)$$
that is the ``space'' of {\it all} $PSL(2,\C)$-valued representations
of the fundamental group $\pi=\pi_1(S\setminus V)$, up to
conjugation. 

\noindent In this section we produce the bundles of cocycle $\Dd$-parameters and
$\Ii\partial$-parameters over $\Rr(g,r)$. It is convenient to replace
$S\setminus V$ with an oriented compact surface with isomorphic
fundamental group. So we fix a compact oriented surface $F$ with $r$
boundary components, obtained by removing from $S$ the interior of a
small $2$-disk $D_i$ such that, for every $i$, $v_i \in \partial D_i$.
Clearly the inclusions of Int$(F)$ into $F$ and $S\setminus V$
respectively induce the identification $\pi = \pi_1(F)$. We stress
that there are no restrictions on the values of the representations at the boundary loops of $F$.

\subsection{Preliminaries on the character variety}

As $\pi$ is a free group with $\kappa = 2g + r -1$ elements, the {\it
  variety of representations} ${\rm Hom}(\pi, PSL(2,\C))$ is naturally
identified with $PSL(2,\C)^\kappa$. Any choice of free generators of
$\pi$ determines such an identification, and the identifications
associated to different choices are related by algebraic automorphisms
of $PSL(2,\C)^\kappa$.  Moreover, the isomorphism $PSL(2,\C)\cong
SO(3,\C)$ implies that ${\rm Hom}(\pi, PSL(2,\C))$ is an {\it
  affine} complex algebraic variety, with the complex algebraic action
of $PSL(2,\C)$. 

\noindent But the quotient space $\Rr(g,r)$ is a much more
delicate object.  This rough topological quotient space is not
even Hausdorff and it is more convenient to consider the algebraic
quotient $X(\pi)={\rm Hom}(\pi, PSL(2,\C))//PSL(2,\C)$ of invariant
theory, called the {\it variety of $PSL(2,\C)$-characters}. We refer
to \cite{HP} for a careful treatment of this matter. We recall that
$X(\pi)$ is an affine complex algebraic set together with a surjective
regular map $t:{\rm Hom}(\pi, PSL(2,\C))\to X(\pi)$, which induces an
isomorphism $t^*$ between the regular functions on $X(\pi)$ and the
regular functions on ${\rm Hom}(\pi, PSL(2,\C))$ invariant by
conjugation. In general $t(\gamma)=t(\sigma)$ does not imply that
$\gamma$ and $\sigma$ are conjugate, i.e. that the quotient set
$\Rr(g,r)$ is $X(\pi)$. However, this is true if we restrict to the
subvariety ${\rm Hom}^{irr}(\pi, PSL(2,\C))$ of ${\rm Hom}(\pi,
PSL(2,\C))$ made by the {\it irreducible} representations (i.e.
without any fixed point in $\mc\mathbb{P}^1$): we have ${\rm
  Hom}^{irr}(\pi, PSL(2,\C)) = t^{-1}(X^{irr}(\pi))$, where
$X^{irr}(\pi)= t({\rm Hom}^{irr}(\pi, PSL(2,\C))$, and the (restricted)
rough quotient $\Rr(g,r)^{irr}$ and the algebraic quotient
$X^{irr}(\pi)$ do agree.

We will not really use the variety of characters, as it is
more convenient for us to consider the variety of representations, by
taking track of the conjugacy action of $PSL(2,\C)$.  More precisely,
we will use algebraic varieties associated to certain {\it
  combinatorial markings} of $F$, that can be considered as
counterparts of the variety of representations, up to suitable {\it
  gauge transformations}, that are counterparts of the conjugacy
action. It shall be useful to consider these objects as
geometric bundles over $\Rr(g,r)$, to treat the ``complex dimension''
and so on. We will do it somewhat formally,
being aware that everything can be substantiated in terms of the
variety of characters, or by restriction to the irreducible
representations. We prefer to treat the whole $\Rr(g,r)$ anyway,
because the construction of the QHFTs does not really require any
restriction on the flat $PSL(2,\C)$-connections. So, as the group
$PSL(2,\C)$ has trivial centre and complex dimension equal to $3$, we
can say that the complex dimension of $\Rr(g,r)$ is equal to $3\kappa
- 3 = -3\chi(F)$. 

\subsection {Efficient triangulations}\label{e-triang}
The first step is to select a class of {\it efficient triangulations}
of $F$. 
\smallskip

We use possibly {\it singular} triangulations of compact oriented
$n$-manifolds $Y$. Any such a triangulation $T$ can be described as a
finite family of oriented {\it abstract} $n$-simplices, together with
the identification of some pairs of abstract $(n-1)$-faces, in such a
way that $Y$ is the quotient space. The face identifications are
orientation reversing so that the orientations of the $n$-simplices
match to produce the given orientation of $Y$. Multiply adjacent as
well as self-adjacent $n$-simplices are allowed.  

Let us start with any {\it branched} triangulation $(T',b')$ of $S$
having $V$ as set of vertices. A {\it branching} $b'$ is a system of
orientations of the edges of $T'$ such that the induced orientations
on the edges of each abstract triangle is compatible with a total
ordering of its vertices, via the rule: each edge is directed towards
the biggest end-point. Hence no abstract triangle of $T'$ inherits
from $b'$ an
orientation of its boundary: only two edges have a compatible {\it
  prevailing} orientation.  If $x_0,x_1,x_2$ are the $b'$-ordered
vertices of a triangle, we name and order its $b'$-oriented edges as:
$e_0=[x_0,x_1], e_1 = [x_1,x_2], e_2 = [x_0,x_2]$, so that $e_0,\ e_1$
have the prevailing orientation.  This also induces a $b'$-orientation
on every triangle, that is the orientation which induces the
prevailing edge orientation. The $b'$-orientation may or may not agree
with the given orientation of $S$. We encode it via a {\it sign
  function} $\sigma = \sigma_{(T',b')}$, defined on the set of
triangles of $T'$, by stipulating that the sign of a triangle is $\pm
1$ if the two orientations do or do not agree respectively. Note that
such a $(T',b')$ exists due to the assumption we have made on the
pair $(g,r)$.

Given any $(T',b')$ as above, we consider {\it corner maps} $v\mapsto c_v$
which associate to each vertex of $T'$ one corner in its star ${\rm
Star}(v)$, and we denote by $v\mapsto t_v$ the induced map that associates
to $v$ the (abstract) triangle that contains the corner $c_v$.  We say
that $v\mapsto c_v$ is {\it $t$-injective} if $v\mapsto t_v$ is injective.
\begin{lem}\label{injcorner} 
For every $(g,r)$ as above, let us assume furthermore that $r>3$ if
$g=0$. Then every triangulation $T'$ of $S$ with $r$ vertices admits
$t$-injective corner maps.
\end{lem}
\noindent
{\it Proof.} First we show that for every $(g,r)$ as in the statement
of the lemma {\it there exist} triangulations of $S$ with $r$ vertices
admitting $t$-injective corner maps. We do it by induction on $r$.
For $(g=0, r=4)$ and $(g>0,r=1)$, it is evident that $T'$ with
$t$-injective $c_v$ do exist.  Clearly, a $t$-injective $c_v$ exists on
any $T''$ obtained from $T'$ via a $1\to 3$ move, i.e. a move that
subdivides one triangle of $T'$ by $3$ triangles, introducing one new
vertex. So we conclude by induction on $r$.

\begin{figure}[ht]
\begin{center}
\includegraphics[width=10cm]{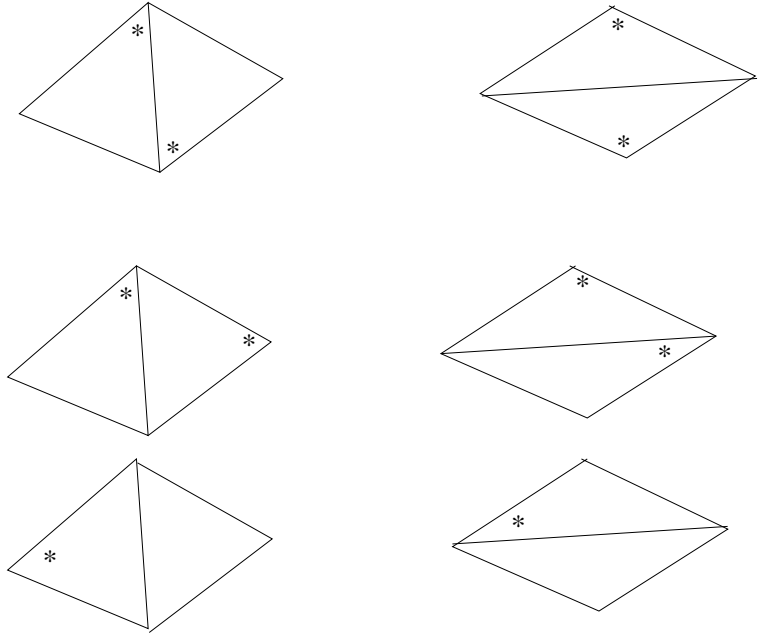}
\caption{\label{eTflip} The flips with marked corners.} 
\end{center}
\end{figure}

\noindent 
In figure \ref{eTflip} the corner selection $v\mapsto c_v$ is specified by
a $*$, and the rows show essentially all possible {\it flips}, up to
some evident variations, that preserve the property that $v\mapsto t_v$ is
injective. Consider any triangulation $T'$ of $S$ with $r$ vertices,
and let $T''$ be a triangulation with the same number of vertices and
which admits a $t$-injective corner map. It is well known that $T''$
is connected to $T'$ via a finite sequence of (naked) flips. The
$t$-injective corner map $v\mapsto c''_v$ for $T''$ transits to a
$t$-injective $v\mapsto c_v$ for $T'$, by decorating these flips as in
Fig. \ref{eTflip}.\hfill$\Box$ \medskip

\noindent The only case excluded by the above lemma is $(g=0,r=3)$;
in this case we have $2$ triangles, hence $t$-injective maps cannot
exists. 

In the generic cases when the lemma applies, let us fix a
$t$-injective corner map $v\mapsto c_v$ for $(T',b')$.  In the
interior of every triangle $t=t_v$ of $T'$ that contains a selected
corner $c_v$ corresponding to a vertex $v$, consider two nested bigons
$ D_v \subset D'_v$ with one common vertex at $v$. Call $v'\in D_v$
and $v''\in D'_v$ the other two vertices of the bigons. Remove from
$t_v$ the interior of $D_v$, obtaining $s_v$. Triangulate $s_v$ by
making the cone with base $v''$. We find a triangulation of $s_v$ with
$5$ triangles, $5$ vertices and $10$ edges. Repeating this procedure
independently on every $s_v$, we get a triangulation $T$ of $F$, with
$3r$ vertices and $p+4r$ triangles, where $p$ denotes the number of
triangles of $T'$.  The set of edges of $T$, $E(T)$, contains $E(T')$,
and $|E(T)|= |E(T')|+7r$.

\begin{figure}[ht]
\begin{center}
\includegraphics[width=10cm]{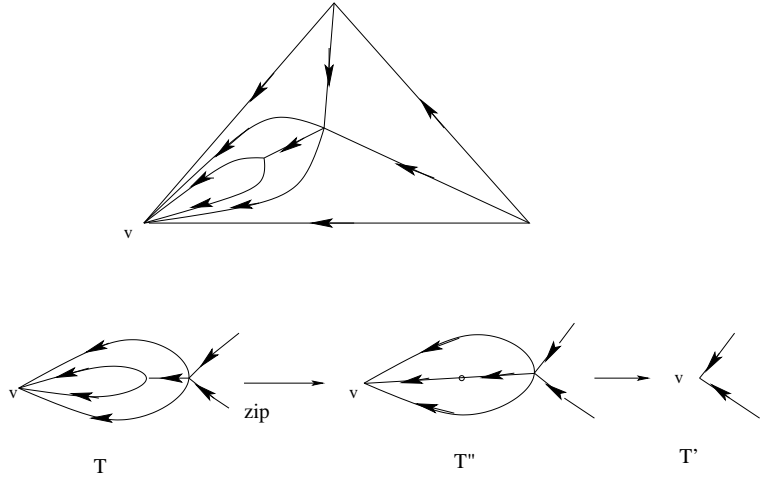}
\caption{\label{newb} The branching $b$.} 
\end{center}
\end{figure}

Now we fix a way of extending the branching $b'$ to a branching $b$ on
$T$. This is shown in in Fig. \ref{newb}. With this choice there is a
clear transition from $(T,b)$ to $(T',b')$: first {\it zip} the two
boundary components of the inner bigon, and get a branched
triangulation $(T'',b'')$ of $S$ with $\bar{V}=V\cup V' \cup V''$ as
set of vertices. Then {\it collapse} each bigon pattern of $T''$ to
the corresponding $v$, and get back the initial $(T',b')$ of $(S,V)$.

\noindent In the generic case, the triangulations $(T,b)$ of $F$
thus obtained are by definition our efficient {$e$-{\it
triangulations}. The triangle sign function $\sigma$ naturally extends
to $(T,b)$.  For each $s_v$, we select a {\it base} triangle among the
$2$ containing a boundary edge. For example, in Fig. \ref{newb}, we
take the triangle $\tau_v$ which contains the boundary edge such that
the $b$-orientation and the boundary orientation do agree. In general we
stipulate that $\sigma (\tau_v)=1$.

\begin{figure}[ht]
\begin{center}
\includegraphics[width=10cm]{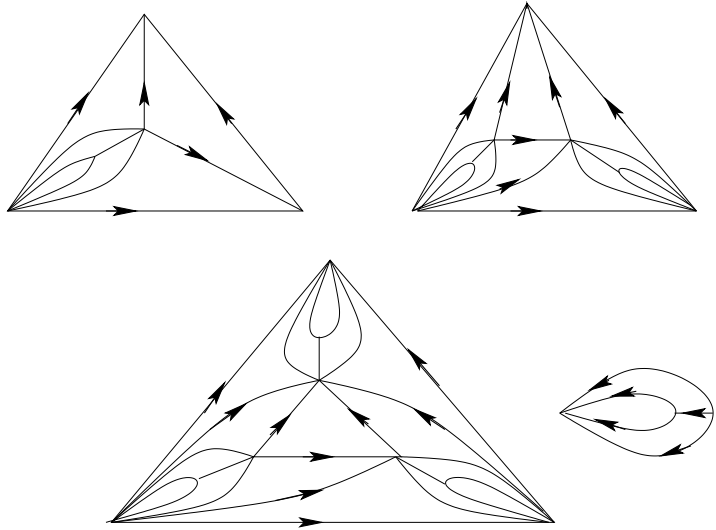}
\caption{\label{Tb} The triangulation $(T,b)$} 
\end{center}
\end{figure}

\noindent 
In the special case $(g=0,r=3)$ we have to consider the further
situation of a triangle containing two selected corners.  This is
shown on the left of the first row of Fig. \ref{Tb}.  In fact this
figure shows essentially {\it all} the possible configurations that we
obtain by using {\it arbitrary} corner maps. The are no conceptual
obstructions to use arbitrary corners maps in what follows. We prefer
to specialize the corner maps just to limit the configurations and
simplify the exposition. Moreover, we will limit ourselves to give the details
in the generic case (and referring to Fig. \ref{newb}), as the
extension to the special case or to the other positions of the
selected corner shall be straightforward.

\subsection{$\Dd$-parameters}\label{Dpar}
Fix an $e$-triangulation $(T,b)$ of $F$. Let us denote by $Z(T,b)$ the
space of $PSL(2,\C)$-valued $1$-cocycles on $(T,b)$. We use
the $b$-orientation of the edges, so that on each triangle with ordered
$b$-oriented edges $e_0,e_1,e_2$ the cocycle condition reads:
$z(e_0)z(e_1)z(e_2)^{-1} = 1$.

We write $C(T,b)$ for the space of $PSL(2,\C)$-valued $0$-cochains, that is
the $PSL(2,\C)$-valued functions defined on the set of vertices of
$T$. Two $1$-cocycles $z$ and $z'$ are said equivalent {\it up to
gauge transformation} if there is a $0$-cochain $\lambda$ such that,
for every (abstract) oriented edge $e=[x_0,x_1]$, we have $z'(e) =
\lambda(x_0)^{-1}z(e)\lambda(x_1)$.  Possibly
$\lambda(x_0)=\lambda(x_1)$, when the two abstract vertices are
identified to one vertex of $T$.

We denote by $H(T,b) = Z(T,b)/C(T,b)$ the quotient set. It is well
known that $H(T,b)$ is in one-to-one correspondence with $\Rr(g,r)$.
More precisely, fix a vertex $x_0$ of $T$ as base point and set $\pi =
\pi_1(F,x_0)$. Then there is a natural surjective map $f: Z(T,b) \to
{\rm Hom}(\pi, PSL(2,\C))$; $f(z)$ and $f(z')$ represent the
same point in $\Rr(g,r)$ iff they are related by gauge
transformations. The complex dimension of $C(T,b)$ is equal to
$3r{\rm dim}(PSL(2,\C))= 9r$. The algebraic set $Z(T,b)$ is defined by
$(p+4r){\rm dim}(PSL(2,\C))= 3(p+4r)$ relations on $|E(T)|{\rm
  dim}(PSL(2,\C)=3|E(T)|$ variables ($p$ is the number of triangles of
the initial triangulation $T'$ of $S$ with $r$ vertices). Remind that
$PSL(2,\C)$ has trivial centre. Hence we find that the (formal)
complex dimension of $H(T,b)$ is just $$3(|E(T)| - (p+4r) - 3r) =
-3\chi(F)$$
that is the dimension of $\Rr(g,r)$. This
essentially means that there are {\it no negligiable} relations
defining $Z(T,b)$.  \smallskip

Let us denote by $p: Z(T,b) \to H(T,b)$ the natural projection. A
way to get {\it $\Dd$-parameters} for $\Rr(g,r)\cong H(T,b)$,
that is parameters based on $1$-cocycle coefficients, should be to
construct nicely parametrized global sections of $p$. Although this is
too optimistic, we will specialize anyhow the cocycles to reduce as
much as possible the set of {\it residual gauge transformations}.

The conjugacy class of every element $g\in PSL(2,\C)$ can be specified by a
symbol $c(g)$, as follows. Set $c(id) := I$, otherwise set either
$c(g)= (a,{\rm diag})$, or
$c(g)= 1$, where: $a \in \C\setminus \{0,1\}$, ``diag'' means that
$g$ is represented (up to conjugation) by a diagonal matrix with $a$
as first eigenvalue; $1$ means that $g$ is represented by the unipotent
upper triangular matrix with $1$ as upper triangular coefficient. In
other words, $c(g)$ determines one distinguished representative in the
conjugacy class of $g$. Sometimes we will say that $g$ is of {\it
trivial}, {\it parabolic} or {\it generic type}, respectively.

In what follows, we will confuse any element of $PSL(2,\C)$ with its
$SL(2,\C)$-representatives. If $B=B(2,\C)$ denote the Borel subgroup
of upper-triangular matrices of $SL(2,\C)$, every $g\in B$ is written
in the form $g=[a,b]$, where $a\in C^*$ is the first eigenvalue of $g$,
and $b$ is the upper-diagonal entry of $g$.

Define
$$ \beta': \Rr(g,r) \to [\{I,1\}\cup (\C^* \times \{{\rm diag}
\})]^r$$ as the map which associates to every holonomy $\rho$ the $r$-uple
$(c(\rho(\gamma_1)),\dots, \ c(\rho(\gamma_r)))$, where $\gamma_i$ is
the oriented boundary loop of $F$ at the vertex $v_i$. Consider
$H(T,b)$ as a set realization of $\Rr(g,r)$, and lift $\beta'$ to $Z(T,b)$
via the composition $\beta=\beta' \circ p$.

For every $x\in \beta(Z(T,b))$, put $Z(T,b,x):= \beta^{-1}(x)$. This is
mapped by $p$ onto $H(T,b,x)=\Rr(g,r,x) = \beta'^{-1}(x)$. So, by
varying $x$, we get a {\it partition} of the projection $p$,
and we consider each piece $p_x: Z(T,b,x)\to \Rr(g,r,x)$.

\smallskip

Denote by $Z(T,b,x)_{\beta}$ the subset of $Z(T,b,x)$ made by the
cocycles $z$ such that, along every oriented boundary loop $\gamma_i$,
the product of the cocycle values, starting from the vertex $v_i$,
exactly equals the distinguished representative $c(p(z)(\gamma_i))$.
Any cocycle $z\in Z(T,b,x)$ can be modified to a cocycle in
$Z(T,b,x)_{\beta}$ via a gauge transformation associated to a suitable
$0$-cochain with support at $V$. So the algebraic set
$Z(T,b,x)_{\beta}$ is non empty and the restriction of $p_x$ maps it
onto $\Rr(g,r,x)$.

\smallskip

\noindent Hence we further restrict ourselves to $p_x: Z(T,b,x)_{\beta} \to
\Rr(g,r,x)$. The set of residual gauge transformations is already
smaller. In fact, to stay in $Z(T,b,x)_{\beta}$ we must act with
$0$-cochains $\lambda$ such that, for every $v_i$,
$\lambda(v_i)$ belongs to the stabilizer ${\rm Stab}(x_i)$ (while
$\lambda$ is arbitrary at the other vertices $v'_i$ and $v_i''$ of
$T$).

By using the notations introduced in Subsection \ref{e-triang}, let us
consider for every $v$ the base triangle $\tau_v$ of $s_v$ with its
$b$-ordered edges $e_0(v), e_1(v)$. Define
$$ \delta_x: Z(T,b,x)_{\beta} \to (PSL(2,\C)^2)^r$$
$$ \delta_x(z) = ([z(e_0(v_1)),z(e_1(v_1))],\dots,\   
[z(e_0(v_r)),z(e_1(v_r))]) \ .$$
For every $y \in (PSL(2,\C)^2)^r$, set $Z(T,b,x,y)_{\beta}:=
\delta_x^{-1}(y)$. As above, any cocycle $z \in Z(T,b,x)_{\beta}$ can
be modified to one in $Z(T,b,x,y)_{\beta}$ just by acting with
a suitable $0$-cochain with support at the vertices $v'_i$ and
$v_i''$'. Hence we have:

\begin{lem} The map $\delta_x$ and the restriction of $p_x$ to $Z(T,b,x,y)_{\beta}$ are surjective.
\end{lem}
We restrict once more to $p_{x,y}: Z(T,b,x,y)_{\beta}\to \Rr(g,r,x)$.
Let us determine the residual gauge transformations. Fix a cell $s_v$
as defined in Subsection \ref{e-triang}.  Consider a $0$-cochain $a$
with support at $v$, $v'$ and $v''$. Set $g_h = z(e_h(v))$, $h=0,1$.
We stay in $Z(T,b,x,y)_{\beta}$ iff:
$$ a(v)\in {\rm Stab}(x_v), \ \ a(v'')^{-1}g_0a(v') = g_0,\ \
a(v')^{-1}g_1a(v) = g_1 \ .$$ Then it is clear that, once $a(v)$ is
fixed, then the rest of the cochain is uniquely determined. As the
$s_v$'s contribute independently each to the other, the set $\Gg
(T,b,x,y)$ of residual gauge transformations is parametrized by
$$\Gg (T,b,x,y)\cong {\rm Stab}(x_1)\times \dots \times {\rm
  Stab}(x_r)\ .$$
Note that we have ${\rm Stab}(I)= PSL(2,\C)$, ${\rm
  Stab}(1) = {\rm Par}(2,\C)$ and ${\rm Stab}(a,{\rm diag}) = {\rm
  Diag}(2,\C)$, the image in $PSL(2,\C)$ of the upper triangular
parabolic (resp. diagonal) subgroup of $SL(2,\C)$. Hence $\Gg
(T,b,x,y)$, in particular its dimension, can be easily determined, and
depend only on the {\it types}, say $t(x)$, of the boundary loops.
\smallskip

\noindent From this we can derive a rather neat qualitative description of
$\Rr(g,r)$. Denote by $Z(T,b)_\beta$ the union of all
$Z(T,b,x,y)_\beta$'s, with the natural projection $p_\beta :
Z(T,b)_\beta \to \Rr(g,r)$. Let
$$
t: [\{I,1\}\cup (\C^* \times \{{\rm diag}\}]^r\to [\{I, 1, \ {\rm
  diag}\}]^r$$
be the natural forgetting map which associates to each
boundary conjugagy class its type.  Define $\phi' = t\circ \beta'$ and
denote by $\phi_\beta$ the restriction to $Z(T,b)_\beta$ of $\phi =
t\circ \beta$.  For every $w\in [\{I,1,\ {\rm diag} \}]^r$ set
$\Rr(g,r,w-type) = \phi'^{-1}(w)$ and $Z(T,b,w-type)_\beta=
\phi_\beta^{-1}(w)$. The above constructions eventually give:

\begin{prop}\label{topcomp} (1) By varying $w$ we get a partition of the projection $p_\beta$ by 
  the maps $p_{w,\beta} : Z(T,b,w-type)_\beta \to \Rr(g,r,w-type)$.
  Moreover, each space $\Rr(g,r,w-type)$ (resp. $Z(T,b,w-type)_\beta$) is
  fibred by the $\Rr(g,r,x)$'s (resp. $Z(T,b,x,y)_\beta$'s) such that
  $t(x)=w$.  \smallskip

(2) For every $y$, the set $Z(T,b,x,y)_\beta$ is a complex affine algebraic
set, which projects via $p_{x,y}$ onto $\Rr(g,r,x)$. We
    have:
$$\begin{array}{lll}
{\rm dim}(Z(T,b,x,y)_\beta)= -3\chi(F)\\ 
{\rm dim}(\Rr(g,r,x)) = -3\chi(F) - {\rm dim}(\Gg (T,b,t(x)))\\
{\rm dim}(\Rr(g,r,w-type)) = {\rm dim}(\Rr(g,r,x)) + \alpha(w)
\end{array}$$
where $w=t(x)$ and $\alpha(w)$ is the number of entries
of $w$ of type ``diag''.
\end{prop}

Let us look at some particular cases:
\smallskip

(a) If $w=w_{{\rm diag}}:=({\rm diag},\dots,\ {\rm diag})$, then ${\rm
  dim}(\Rr(g,r,w_{{\rm diag}}-type))= -3\chi(F)$. 

\smallskip

(b) If $w=w_{1}:=(1,\ \dots,\ 1)$, then ${\rm dim}(\Rr(g,r,w_1-type))= -3\chi(F)- r$.  

\smallskip

(c) If $w=w_I:=(I,\ I,\ \dots,\ I)$, then ${\rm dim}(\Rr(g,r,w_I-type))= -3\chi(F) -3r = -3\chi(S) = 6g-6$.

\medskip

We say that: a type $w'$ is obtained from $w$ via a {\it simple
  degeneration} $w\to w'$ if they differ at just one entry, where
either ${\rm diag}\to 1$ or $1 \to I$; $\Rr(g,r,w'-type)$ is {\it in
  the formal frontier} of $\Rr(g,r,w-type)$ if $w'$ is obtained from
$w$ by a finite sequence of simple degeneration. In such a case,
clearly ${\rm dim}(\Rr(g,r,w'-type))<{\rm dim}(\Rr(g,r,w-type))$.
Hence (a) above implies that $\Rr(g,r,w_{{\rm diag}}-type)$ is a dense open set in
$\Rr(g,r)$, and we have a nice {\it filtration} of our phase space
for which $\Rr(g,r,w_I-type)$ is the `deepest' part. Presumably this
filtration can be refined to a nice {\it stratification}, via the
study of the singularities of each $\Rr(g,r,w-type)$ and of its actual
closure in $\Rr(g,r)$.

\smallskip

Summarizing, it is convenient to work separately on each
$\Rr(g,r,w-type)$, more precisely on each $\Rr(g,r,x)$ such that
$t(x)=w$. The idea is that, for a fixed $w$, we get one ``bundle''
$Z(T,b,w-type)_\beta \rightarrow \Rr(g,r,w-type)$ {\it fibred over the
  parameters $x$}. To concretize this idea, we will
suitably specify in the next discussion about $\Ii$-parameters a
choice $y_x$ for $y$, depending on $x$, so that we get:

\begin{defi}\label{cocpar}{\rm Set $\bar{p}_x=p_{x,y_x}$. The bundle
    of {\it $w$-type
    cocycle $\Dd$-parameters} is given by}
$$\coprod_{x \in t^{-1}(w)} \bar{p}_x :
    Z(T,b,x,y_x)_{\beta}\rightarrow \Rr(g,r,w-type).$$
\end{defi}
The equivalence
classes of cocycle $\Dd$-parameters up to residual gauge
transformations give coordinates for $\Rr(g,r,x)$.

\subsection {$\Ii$-parameters}\label{Ipar}
Let $(T,b)$ be an $e$-triangulation of $F$ as above. Recall that:
$(T,b)$ is obtained from a branched triangulation $(T',b')$ of $S$; 
$\chi (F)=\chi(S)-r$; $T'$ has $r$ vertices and $p$ triangles, so that 
$2|E(T')|=3p = -6\chi(F)$. The triangulation $T$ has $3r$ vertices and $p+4r = 4r-2\chi(F)$
triangles, and $|E(T)|=7r - 3\chi(F)$.  
\medskip

Before proceeding, we recall few facts about $3$-{\it dimensional
branchings} (see \cite {BB1, BB2} for full details). 

\paragraph {Branchings.}\label{3Dbranch}
Given a triangulation $T$ of an oriented compact $3$-manifold $Y$, a
branching $b$ on $T$ is a system of orientations of the edges of $T$
which induces on each abstract tetrahedron $\Delta$ of $T$ a total
ordering $x_0,x_1,x_2,x_3$ of its vertices. Note that each $2$-face
of $\Delta$ inherits a $2$-dimensional branching in the sense just
defined. The ambient orientation of $Y$ induces an orientation on each
$\Delta$. Also the branching induces a $b$-orientation on $\Delta$:
the $b$-orientation coincides with the ambient orientation iff the
$b$-orientation of the $2$-face $f(3)$ opposite to the vertex $x_3$ of
$\Delta$ coincides with the boundary orientation (i.e.
$\sigma(f(3))=1$). The {\it sign function} $*_b$ for the tetrahedra of
$(T,b)$ is defined by $*_b(\Delta) := \sigma (f(3))$. For every
$(\Delta,b)$ we denote by $e_0,e_1,e_2$ the $b$-ordered and oriented
edges of $f(3)$.

\noindent
It is convenient to give also an encoding of these $3$-dimensional
branched triangulations $(T,b)$ in terms of their {\it dual cell
  decompositions}. In Fig. \ref{2CQDidealtensor} we see an enriched
version of the $1$-skeleton of such a dual cell decomposition,
localized at branched tetrahedra of sign $*_b=\pm 1$ (ignore the
symbols $x$, $\alpha, \dots, \delta$ for the moment, as they refer to
later considerations).

\begin{figure}[ht]
\begin{center}
\includegraphics[width=10cm]{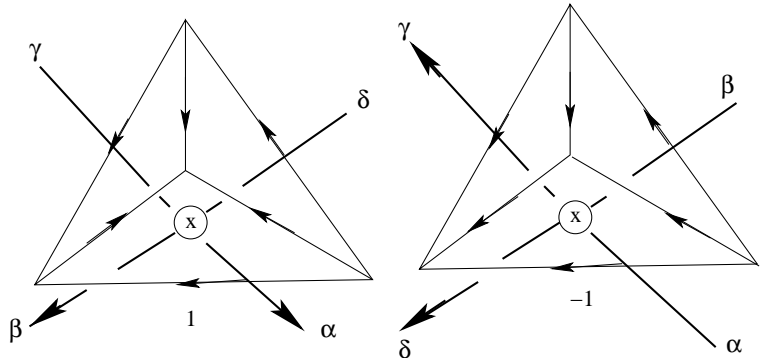}
\caption{\label{2CQDidealtensor} Decorated tetrahedra and dual
encoding.}
\end{center}
\end{figure}

In the picture we see in fact a {\it planar} realization of this
$1$-skeleton; its four branches at the vertex dual to int$(\Delta)$
are arranged to form a normal crossing with an under/over arc
specification (like for ordinary link diagrams). The $b$-sign
$*_b=\pm1$ is encoded by the usual normal crossing index. We have
omitted to draw any arrows on two branches, as it is inderstood that they
are incoming at the crossing. Note that these decorated graphs contain
all the information in order to reconstruct the corresponding dual
branched tetrahedron $(\Delta,b)$. The oriented branches are outgoing exactly when
the corresponding dual $2$-face has $b$-sign $\sigma = 1$.

\paragraph {The cylinder $\Cc(T,b)$.}\label{cylinder}
Consider the cylinder $\Cc=\Cc_F= F\times [-1,1]$, oriented in such a
way that the oriented surface $F$ is identified with the oriented
``horizontal'' boundary component $F_ - = F\times \{-1\}$ of $\Cc$.
Sometimes we write $(T_-,b_-)$ for the branched triangulation on this
boundary component, via this identification. Similarly we
write $(T_+,b_+)$ for the copy of $(T,b)$ on the other horizontal boundary
component $F_+$, which has the opposite orientation. Finally, we denote by
$\Cc(T,b):=(\Cc(T),\Cc(b))$ the branched triangulation of $\Cc$
obtained as follows. \smallskip

Consider first the natural {\it product} cell decomposition $P(T,b)$
of $\Cc$, made by $p+4r$ prisms with triangular base. We stipulate that
all the $3r$ vertical edges of $P(T,b)$ are oriented to point towards
$F_+$. For every abstract prism $P$, every ``vertical'' quadrilateral,
$R$ say, on its boundary has both the two horizontal and the two
vertical edges endowed with parallel orientations. So exactly one
vertex of $R$ is a source (that belongs to $F_-$), and exactly one is a
pit (that belongs to $F_+$). Then we triangulate each $R$ via the
oriented diagonal going from the source to the pit. Finally we extend
the so obtained triangulation of $\partial P$, to a triangulation of
$P$ by $3$ tetrahedra, just by making the cone from the $b$-first
vertex of the bottom base triangle of $P$ (note that no further
vertices nor further edges have been introduced). Repeating this for
every prism, we finally get our branched triangulation $\Cc(T,b)$ of $\Cc$. 

\smallskip

Let us list few properties of $\Cc(T,b)$:

\smallskip

(1) $\Cc(T,b)$ contains $6r$ vertices (that are all on $F_{\pm}$) and
$3(4r-2\chi(F))$ tetrahedra; moreover, it has $2(7r - 3\chi(F))$
{\it horizontal} edges on $F_{\pm}$, $r$ {\it vertical} edges over
the vertices of $T$, and $7r - 3\chi(F)$ {\it diagonal} edges on the
vertical rectangles. Note that each prism contains $10$ triangles.
\smallskip

(2) For every prism $P$, denote:

- $t_{\pm}= t_{\pm}(P)$ the base triangle contained in $F_{\pm}$;

- $\Delta_{\pm}= \Delta_{\pm}(P)$ the tetrahedron based at $t_{\pm}$;

- $\Delta_0 = \Delta_0(P)$ the interior tetrahedron.
\smallskip

\noindent Assume that $\sigma_b(t_-)=1$ (so that $\sigma_b(t_+)= -1$).
Then both $*_{\Cc(b)}(\Delta_{\pm})=1$. The tetrahedron $\Delta_0$
shares one edge with each base triangle respectively; these are
opposite edges of $\Delta_0$. We have $*_{\Cc(b)}(\Delta_0)=-1$. The
$b$-oriented dual graph of $\Cc(T,b)|P$ points outside at $t_-$, and
points inside at $t_+$.  If $*_b(t_-)=-1$ the same facts hold,
providing that all the signs are inverted.

\smallskip

The $\Ii$-parameters for $\Rr(g,r)$ shall result from a suitable {\it
idealization} procedure of the $\Dd$-parameters discussed in
Subsection \ref{Dpar}. As explained there, we will work
separately on each $\Rr(g,r,x)$ by using the surjective
projections $p_{x,y}: Z(T,b,x,y)_{\beta}\to \Rr(g,r,x)$. The rough
idea is to extend each $z\in Z(T,b,x,y)_{\beta}$ to some cocycle
$\Cc(z)\in Z(\Cc(T,b))$ and take (if possible) its idealization,
thus obtaining the corresponding cross-ratio $\Ii$-parameters.

\paragraph {Idealization.}\label{idealization}
Let us briefly recall few general facts about the idealization
procedure.  We refer to \cite{BB1, BB2} for the full details.  Let
$(T,b,z)$ be a branched triangulation of an oriented compact
$3$-manifold $Y$, equipped with a $PSL(2,\C)$-valued $1$-cocycle $z\in
Z(T,b)$.  We fix once for ever $0\in \C$ as {\it base point} of our
idealization procedure. We say that an abstract tetrahedron
$(\Delta,b,z)$ of $(T,b,z)$ (with the induced branching and cocycle), is {\it idealizable} if
$$
u_0=0,\ u_1= z_0(0),\ u_2= z_0z_1(0),\ u_3= z_0z_1z'_0(0) $$
are
$4$ distinct points in $\mc \subset \mc\mathbb{P}^1= \partial \bar{\mh}^3$. Here
$z_i = z(e_i)$. These $4$ points span a (possibly flat) {\it
  hyperbolic ideal tetrahedron with ordered vertices}.  We call
$(T,b,z)$ a $\Dd$-{triangulation} if all its tetrahedra are
idealizable. The {\it idealization} $(T,b,w)$ of a $\Dd$-triangulation
consists of the family $\{(\Delta,b,w)\}$, where $\Delta$ spans the
$3$-simplices of $T$, and each edge $e$ of $\Delta$ is now decorated by
the the appropriate cross-ratio modulus $w(e)\in \C \setminus
\{0,1\}$ of the above hyperbolic ideal tetrahedron. In fact the $w(e)$'s are specified by
the {\it modular triple} $w=(w_0,w_1,w_2)$, $w_i = w(e_i)$, as
opposite edges share the same cross-ratio modulus. The idealization
$(T,b,w)$ gives a so called $\Ii$-{triangulation} of $Y$. This means that at each {\it
  internal} (i.e. not contained in $\partial Y$) edge $e$ of $T$ it is
satisfied the {\it edge compatibility condition}
\begin{equation}\label{ideq1defsec}
\prod_{h \in \epsilon_T^{-1}(e)} w^j(h)^{*_{b^j}} =1
\end{equation}
where $\epsilon_T$ is the map that
associates to every abstract edge the corresponding edge in $T$ (via
the face identifications), and $*_{b^j}=\pm 1$ according to the
$b^j$-orientation of the tetrahedron $\Delta^j$ that contains the
abstract edge $h$.

\smallskip

$\Ii$-triangulations actually encode $PSL(2,\C)$-valued
representations of the fundamental group of $Y$ up to conjugation.
More precisely, by lifting a given $\Ii$-triangulation of $Y$ to its
universal covering $Y'$, we can construct, by ``developing'' the
hyperbolic ideal tetrahedra of the triangulation in the naturally
compactified hyperbolic space $\bar{\mh}^3$, a {\it pseudo developing
  map} $d: Y' \to \bar{\mh}^3$ and a representation $h: \pi_1(Y) \to
PSL(2,\C)$, such that $d(\gamma(y)) = h(\gamma)(d(y))$ for every
$\gamma \in \pi_1(Y)$ and $y\in Y'$. The pseudo-developing map $d$ is
unique up to post-composition with the action of $PSL(2,\C)$ on
$\bar{\mh}^3$, and $h$ is unique up to conjugation. 
\begin{defi} \label{triideal} {\rm We say that a cocycle $z \in Z(T,b,x,y)_{\beta}$ is
{\it idealizable} if for any triangle
    of $(T,b)$ the points $u_0=0,\ u_1= z_0(0)$ and $u_2= z_0z_1(0)$ are
    distinct in $\partial \bar{\mh}^3$. We denote by
    $Z_I(T,b,x,y)_{\beta}$ the set of these idealizable cocycles.}
\end{defi}
Note that if there exists an idealizable extension $\Cc(z)$ of $z$ to $C(T)$,
then $z\in Z_I(T,b,x,y)_{\beta}$. We construct such extensions as
follows. 

\smallskip

- Take first the following {\it trivial extension} $\Cc^0(z)$
of $z$ to $\Cc(T,b)$.

- Copy $z$ on the triangulations $(T_{\pm},b_{\pm})$ 
of the horizontal boundary components $F_{\pm}$ of $\Cc$.

- Every vertical quadrilateral $R$ of $\Cc(T,b)$ has the bottom
and top horizontal edges endowed with parallel orientations and with
the same cocycle value, say $g$. Give each vertical edge the value
$1=[1,0]$; there is a unique way to complete the cocycle, just by
giving each diagonal edge the corresponding value $g$.  \smallskip

Evidently $\Cc^0(z)$ is {\it not} idealizable. We have to perturb
it. For every $a\in \mc^*$, consider the $0$-cochain $c^-_a$ that gives
each vertex of $T_-$ the value $[1,0]$, and each vertex of $T_+$ the
value $[1,a]$. Finally let $\Cc^-(z,a)$ be the cocycle on $\Cc(T,b)$
obtained by perturbing $\Cc^0(z)$ via the gauge transformation
corresponding to $c^-_a$.

\smallskip

The proof of the following lemma is easy. Recall that a triangulation
is said {\it quasi-regular} if every edge has {\it distinct}
end-points.

\begin{lem} \label{defIco} For every $(T,b)$, $x$, $y$ and $a$ as above we have:
  \smallskip
  
  (1) $Z_I(T,b,x,y)_\beta$ is a non empty dense open subset of
  $Z(T,b,x,y)_\beta$.

(2) If $(T,b)$ is quasi-regular, then the projection of
$Z_I(T,b,x,y)_\beta$ covers the whole of $\Rr(g,r,x)$.

(3) We can remove from $Z_I(T,b,x,y)_\beta$ a finite number of complex
algebraic hypersurfaces, to obtain an algebraic set $Z_{I,a}(T,b,x,y)_\beta$
such that for every $z$ in $Z_{I,a}(T,b,x,y)_\beta$ the cylinder
cocycle $\Cc^-(z,a)$ is idealizable.

(4) There is a finite number of non zero complex numbers $s$ such that the corresponding
$Z_{I,s}(T,b,x,y)_\beta$'s cover the whole of $Z_I(T,b,x,y)_\beta$.
\end{lem}  
We are interested to the portion $\Rr_I(g,r,x)$ of $\Rr(g,r,x)$
covered by the projections of these $Z_I(T,b,x,y)_\beta$'s. In order
to make everything more definite, we are going now to specify a
normalized choice $y_x$ of $y$.
\paragraph {Normalizing $y=y_x$.}\label{normaly}
Let $S$, $F$ (identified with $F_-$), $(T',b')$, with associated
$e$-triangulation $(T,b)$ of $F$, be as usual. For every vertex $v$ of
$(T',b')$, consider the corresponding cell $s_v$ containing the other two
vertices $v'$ and $v''$ of $T$. Recall the base triangle $\tau=\tau_v$
defined in Subsection \ref{e-triang}. Denote by $\tau'=\tau'_v$ the 
triangle in $s_v$ that contains the other boundary
edge.  Call $e_0=e_0(v),e_1=e_1(v),e_2=e_2(v)$ the $b$-ordered edges
of $\tau$. Set similarly $e'_j$ for $\tau'$. Note that $e_0 = e'_0$.
\smallskip

\noindent All $z\in Z(T,b,x,y)_{\beta}$ share, by definition, the
same contribution of $y$ at $v$, that is
$y_v = (g_0,g_1):=(z(e_0),z(e_1))$.  Set $g_2=g_0g_1$.  Similarly
denote $h_0,h_1,h_2$ on $\tau'$. Recall that $x_v=(h_1)^{-1}g_1$, and that
\smallskip

- $x_v=[s,0]$, $s\in \C\setminus\{0,1\}$, in the case of
{\it generic} boundary holonomy;

- $x_v =[1,1]$, in the case of {\it parabolic} boundary holonomy;

- $x_v =[1,0]$, in the {\it trivial} case.
\smallskip

\noindent 
For any $0$-cochain $c$ with support in $s_v$, write $c=c(v)$,
$c'=c(v')$ and $c''=c(v'')$. If $c$ leads to residual gauge
transformations for $Z(T,b,x,y)_{\beta}$, then $ c\in {\rm Stab}(x_v)$, that is 
\smallskip

- $c=[d,0]$, $d\in \C^*$, in the generic case;

- $c=[1,b]$, $b\in \C$, in the parabolic case;

- $c$ is an arbitrary element of $PSL(2,\mc)$ in the trivial case.
\smallskip

\noindent {\it Normalization in the generic case.}  Set $g_0 = [1,1]$
and $g_1 = [f,0]$, with $f\in \C\setminus \{0,1\}$. Then $g_2 =
[f,1/f]$, $c' = g_1cg_1^{-1}= [f,0][d,0][1/f,0] = c$ and
$c''=g_0c'g_0^{-1}= [1,1][d,0][1,-1]= [d,(1-d^2)/d]$.  Moreover, $h_0
= g_0 = [1,1]$, $h_1= g_1[1/s,0]= [f/s,0]$, and $h_2 = h_0h_1 =
[1,1][f/s,0]= [f/s, s/f]$.

\noindent 
Hence the idealization of $\tau$ has vertices $(0,1,f^2)$. The
idealization of $\tau'$ has vertices $(0,1, (f/s)^2)$. To have in both
cases $3$ distinct points we have only to impose that $f\ne \pm s$ and
$f \ne \pm 1$. We get our normalization by setting $f^2= s$ and
taking the determination of the square root associated to the branch
of logarithm with arguments in $]-\pi,\pi]$.  \smallskip

\noindent {\it Normalization in the parabolic and trivial cases.}
Consider the parabolic case. We manage similarly in order to get first
that $c=c'$.  Put $g_0=[1,1]$, $g_1 = [1,f]$ with $f\ne 0,-1$, and $g_2 =
[1,f+1]$. As $c=[1,b]$, then $c'=[1,f][1,b][1,-f]=[1,b]=c$. Also, we
have $h_0=[1,1]$, $h_1=g_1[1,-1]= [1,f-1]$, and $h_2 = [1,f]$. The
vertices of the idealization of $\tau$ are $(0,1,1+f)$, the ones of
$\tau'$ are $(0,1,f)$. In the case of a trivial boundary loop we have
$g_i = h_i$, $i=0,1,2$. Then we can impose that by
working both in the generic and parabolic ``styles'', we eventually
get the same idealization for $\tau$. So we impose $f^2=1+f$, that is
$f=(1+\sqrt{5})/2$. We take the same normalization also in the
parabolic case.  \medskip

For every $x$, we denote by $Z(T,b,x,y_x)_{\beta}$ the subset of
$Z(T,b,x)_{\beta}$ obtained by performing on each $s_v$ the above
normalization. The surjective projections $\bar{p}_x=p_{x,y_x}:
Z(T,b,x,y_x)_{\beta}\to \Rr(g,r,x)$ form the bundle of cocycle
$\Dd$-parameters of Definition \ref{cocpar}. From now
on we will apply the previously discussed idealization procedure to
this normalized situation.

\paragraph{Construction of the $\Ii$-parameters.} The image of the
sets $Z_{I,a}(T,b,x,y_x)_\beta$, defined in Lemma \ref{defIco}, by the projections $\bar{p}_x$ give our favourite {\it patches} for $\Rr_I(g,r,x)$. On each
$Z_{I,a}(T,b,x,y_x)_\beta$ we have a ``change of coordinates''
$$i: Z_{I,a}(T,b,x,y_x)_\beta\rightarrow \Ii_a(T,b,x,y_x)_\beta$$
that passes
from the cocycle $\Dd$-parameters to {\it $\Ii$-parameters}, namely the
cross-ratio moduli of the tetrahedra in the cylinder $\mathcal{C}(T,b)$,
obtained via the idealization of the $\Cc^-(z,a)$'s. Every space of
$\Ii$-parameters $\Ii_a(T,b,x,y_x)_\beta$ is a Zariski open set of an
algebraic subvariety of
$$(\C \setminus \{0,1\})^{3(2r-2\chi(F))}.$$ Indeed, there are in total
$3(4r-2\chi(F))$ tetrahedra in $\mathcal{C}(T,b)$, but the moduli of the tetrahedra over
each pair of triangles at the boundary loops of $F$ are fixed by the
normalization.
 
\noindent 
This subvariety is defined by the following set of algebraic equations
:
\smallskip

- {\it Diagonal relations.} These are $|E(T)|-3r =4r - 3\chi(F)$ edge
compatibility conditions (\ref{ideq1defsec}) at the internal diagonal
edges of the vertical quadrilaterals. Again because of the
normalization, we can ignore, for every vertex $v\in V$, the $3$
quadrilaterals that lie over the two loop boundary edges and over the
edge connecting $v'$ and $v''$ respectively.

\smallskip

- {\it Vertical relations.} These are $2r$ relations at the vertical
edges of $\Cc(T,b)$ over the $v$- and $v''$-vertices respectively. The
second ones are again $\Ii$-triangulation edge compatibility
conditions (\ref{ideq1defsec}) at these interior edges. The first ones
are also edge compatibility conditions, once we have filled each
vertical boundary tube of $\partial F \times [-1,1]$ by a suitable
$\Ii$-cusp, so that also the vertical edges over the $v$-vertices
become interior edges (the needed {\it cusp} machinery is developed in
the next Section \ref{CUSP}). These relations actually depend on the
values $x_v$ at the corresponding boundary loops.

\noindent So, we have in total $6r - 3\chi (F)$ relations defining
$ \Ii_a(T,b,x,y_x)_\beta$. This gives
$$
{\rm dim} (\Ii_a(T,b,x,y_x)_\beta) \geq 3(2r-2\chi(F))- (6r -
3\chi(F)) = -3\chi(F)\ .$$
On the other hand, $\Ii_a(T,b,x,y_x)_\beta$
is an open set of a space of complex dimension $-3\chi(F)$. It
projects onto $\Rr_I(g,r,x)$, since $\Rr(g,r)$ is encoded
by the $\Ii$-triangulations of the cylinder $C_F$. Hence we eventually
get the following remarkable facts:

\begin{prop} \label{Ipar1} (1) We have ${\rm dim}
  (\Ii_a(T,b,x,y_x)_\beta)= -3\chi(F)$, i.e. the system of equations
  defining $\Ii_a(T,b,x,y_x)_\beta$ is not overdetermined.

(2) The residual $\Dd$-gauge transformations $\Gg(T,b,t(x))$ transit
via the idealization map $i: Z_{I,a}(T,b,x,y_x)_\beta\to
\Ii_a(T,b,x,y_x)_\beta$ onto a space of residual $\Ii$-gauge
tranformations $\Gg_I(T,b,t(x))$ of
the same dimension. So we have a {\it bundle of cross-ratio $\Ii$-parameters}
$$p_{I,a}: \Ii_a(T,b,x,y_x)_\beta \rightarrow \Rr_I(g,r,x)$$
with structural group $\Gg_I(T,b,t(x))$.
\end{prop}

\begin{remark}\label{real} {\rm
We can replace $PSL(2,\C)$ by $PSL(2,\R)$, and almost everything
can be repeated verbatim. A main difference is that in the case of
parabolic ends we have $c(g)=\pm 1$, that is to every parabolic end
there is an associated {\it sign}. Moreover, via the idealization,
we get only  degenerate tetrahedra, that is only {\it real} $\Ii$-moduli.}
\end{remark}

\paragraph{The $\Ii$-boundary map.} We are interested now to the two-dimensional ``trace'' on
$(T,b)$ of the cross-ratio $\Ii$-parameters. For every $a\in \C^*$,
we can use the same formula which enters the edge compatibility
condition (\ref{ideq1defsec}) for the interior edges of any
$\Ii$-triangulation, to define a map
$$ W_a : Z_{I,a}(T,b,x,y_x)_\beta \to (\C\setminus \{0,1\})^{|E(T)|-3r}$$
which factorizes via the idealization as
$$W_a=W^I_a \circ i : Z_{I,a}(T,b,x,y_x)_\beta \to
\Ii_a(T,b,x,y_x)_\beta \to (\C\setminus \{0,1\})^{|E(T)|-3r}.$$
More precisely, as
usual let us identify $(T,b)$ with $(T_-,b_-)$, so that the edges of $T$
are contained in the horizontal boundary component $F_-$ of $\Cc$. The
map $W^I_a$ associates to each edge $e$ of $T=T_-$ the signed product
\begin{equation}\label{ideq2defsec}
\prod_{h \in \epsilon_T^{-1}(e)} w^j(h)^{*_{b^j}} 
\end{equation}  
of the cross-ratio moduli of $i(\Cc^-(z,a))$ at the abstract edges $h$
of $C(T)$ descending onto $e$. (We can impose in Definition \ref{triideal} the genericity condition
that $W_a^I(z)(e)\ne 1$ for all $e$). Note that for every $v \in
V$, the products (\ref{ideq2defsec})
for the two boundary edges of $F$ at $v$ and the edge connecting
$v'$ and $v''$ are fixed by the normalization.

\noindent Here is an interpretation of the $W_a(z)(e)$'s. Consider the two
(branched) triangles adjacent to an edge $e$ of $T$. Let us order, if
needed, the two vertices opposite to $e$ by orienting the edge $e'$
connecting them (i.e. dual to $e$), so that $(e',e)$ defines the
positive orientation of $F$. Then the vertices in the star of $e$
are totally ordered. As in the above discussion about the idealization
procedure, consider the $z$-orbit of $0$ along the oriented edges of
both triangles. The meaning of $W_a(z)(e)$ is that of a cross-ratio
for the quadrilateral in $\mc\mathbb{P}^1 =\partial \bar{\mh}^3$, with
ordered vertices the so obtained four orbit points.  We have:

\begin{lem}\label{W} The maps $W_a$'s match on the overlaps of the
  spaces $Z_{I,a}(T,b,x,y_x)_\beta$, so that we have a well defined {\rm
    $\Ii$-boundary map}
$$ W : Z_I(T,b,x,y_x)_\beta \to (\C\setminus \{0,1\})^{|E(T)|-3r}.$$  
\end{lem} 

\noindent {\it Proof.} We have to show that  
two cocycles $\Cc^-(z,a)$ and $\Cc^-(z,a')$ that differ only for the
coefficients $a$ and $a'$ lead to $W^I_a= W^I_{a'}$. Consider the
gluing $\Cc(T,b)\cup -\Cc(T,b)$, where we have inverted the
orientation of the second copy, and the gluing is made by identifying
the two copies of $F_-$. Also $\Cc^-(z,a)$ and $\Cc^-(z,a')$ glue
together and give us an idealizable cocycle on $\Cc(T,b)\cup
-\Cc(T,b)$. Every edge $e$ on $F_-$ is now an interior edge, and the
usual edge compatibility condition exactly means that
$W_a^I(e)W^I_{a'}(e)^{-1} = 1$.\hfill$\Box$

\medskip

Recall that $\Ii_a(T,b,x,y_x)_\beta$ is an open affine algebraic set
of complex dimension $-3\chi(F)$. As the maps $W^I_a$ are given by
monomials, ${\rm Im}(W)$ is also open with ${\rm dim}({\rm Im}(W))\leq
-3\chi(F)$. Moreover, $\Rr_I(g,r,x)$ is encoded just by the products
(\ref{ideq2defsec}), rather than the whole set of cross-ratio
$\Ii$-parameters (see the next paragraph). A direct computation shows
that the group $\Gg(T,b,t(x))$ of residual $\Dd$-gauge transformations
transit via the map $W_a$ to a group of the same dimension. So we eventually get:

\begin{prop} \label{Ibpar} The set ${\rm Im}(W)$ is an affine complex algebraic
  set of dimension $-3\chi(F)$, that is the total space
  of a bundle over $\Rr_I(g,r,x)$.
\end{prop}

\begin{defi} {\it We call $W(T,b,x)={\rm Im}(W)$ the space of {\it
      $\Ii\partial$-parameters} for $\Rr(g,r,x)$, and denote
      $\pi:W(T,b,x) \rightarrow \Rr_I(g,r,x)$ the associated bundle.}
\end{defi} 

Recall from Lemma \ref{defIco} (2) that if $T$ is quasi-regular, then
$\Rr_I(g,r,x)=\Rr(g,r,x)$. We note that the relations between the
$\Ii\partial$-parameters are very implicit compared to the very
transparent edge compatibility relations in $\Ii_a(T,b,x,y_x)_\beta$.

\paragraph{Holonomies from $\Ii\partial$-parameters.} We can describe explicitely
the bundle map $\pi:W(T,b,x) \rightarrow \Rr_I(g,r,x)$. For a
conjugacy class of representations $\rho\in \Rr(g,r,x)$, take an
$e$-triangulation $T=T_-$ for the surface $F=F_-$ (viewed as the lower
horizontal boundary component of $\mathcal{C}$), such that there
exists an idealizable cocycle $z \in Z_I(T,b,x,y_x)_{\beta}$
representing $\rho$. For instance, such a $z$ exists if $T$ is quasi-regular;
if $\rho$ is quasi-Fuchsian, or more generally if there exists a
non-empty domain of $\mc\mathbb{P}^1$ on which $\rho(\pi_1(F))$ acts
freely, then any $T$ works. We construct representatives
$\tilde{\rho}$ of $\rho$ from any point $W=W(z)$ in the fiber
$\pi^{-1}(\rho)$ as follows.

\begin{figure}[ht]
\begin{center}
\includegraphics[width=10cm]{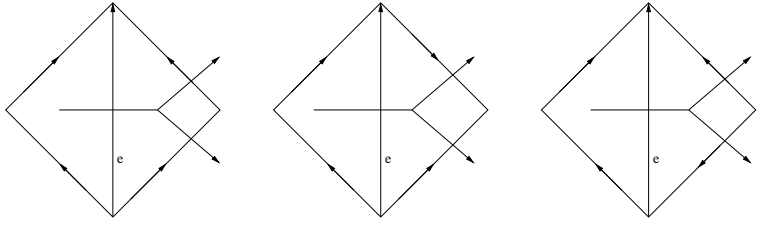}
\caption{\label{matrixrep} The recipe for reading off holonomies from $\Ii\partial$-parameters.}
\end{center}
\end{figure}

Choose a base point $q$ in $F$ not in the $1$-skeleton of $T$. Given
an element of $\pi_1(F,q)$, represent it by a closed curve $\gamma
\subset F$ transverse to $T$, and which do not backtrack (i.e. it never
departs from an edge it just entered). Assume that $\gamma$ intersects
an edge $e$ of $T$ positively w.r.t. the orientation of $F$. Fig.
\ref{matrixrep} shows the three possible branching configurations for
the two triangles glued along $e$.

Fix arbitrarily a square root $W'(e)$ of $W(e)$. Consider the elements
of $PSL(2,\mc)$ given by
$$\gamma(e)=\left( \begin{array}{cc} W'(e) & 0 \\ 0 & W'(e)^{-1}
  \end{array} \right) \quad , \quad p=\left(
  \begin{array}{cc} 0 & 1 \\ 1 & 0  \end{array}\right)\quad ,
\quad l=\left(\begin{array}{cc} -1 & 1 \\ -1 & 0 
\end{array}\right)$$
and $r=l^{-1}$. Recall that $PSL(2,\mc)$ is isomorphic to
Isom$^+(\mh^3)$, with the natural conformal action on
$\mc\mathbb{P}^1=\partial\bar{\mh}^3$ via linear fractional
transformations. The matrix $\gamma(e)$ represents the isometry with
fixed points $0,\infty \in \mc\mathbb{P}^1$ and mapping $1$ to $W(e)$.
The elliptic elements $p$ and $l$ send $(0,1,\infty)$ to
$(\infty,1,0)$ and $(\infty, 0,1)$ respectively. 

The flat principal $PSL(2,\mc)$-bundles associated to $\rho$ carry
parallel transport operators, that we may compute along $\gamma$ by
using the cocycle $z$. For the portion of $\gamma$ represented in the
left configuration of Fig. \ref{matrixrep}, if $\gamma$ turns to the
left after crossing $e$ the parallel transport operator is
$\gamma(e)\cdot p\cdot l$, while it is $\gamma(e)\cdot p\cdot r$ if
$\gamma$ turns to the right. (The composition is on the right, as is
the action of $PSL(2,\mc)$ on the total spaces of the bundles given by
$\rho$).  Similarly, in the middle and right pictures the parallel
transport operators are given by $\gamma(e)\cdot l$ or $\gamma(e)\cdot p\cdot
l$, and $\gamma(e)\cdot p\cdot r$ or $\gamma(e)\cdot r$ respectively.
The action of $p$, $l$ and $r$ depends on the reordering of the
vertices after the mapping $\gamma(e)$. Note that the whole branching
configuration enters the computation of $W(e)$.  If $\gamma$
intersects $e$ negatively, essentially we have to replace $\gamma(e)$
with $\gamma(e)^{-1}$ in the above expressions.  Continuing this way
each time $\gamma$ crosses an edge of $T$ until it comes back to $q$,
we get an element of $PSL(2,\mc)$.

This element does depend only on the homotopy class of $\gamma$ (based
at $q$), for any $W \in \pi^{-1}(\rho)$, and we eventually obtain a
well-defined representation $\tilde{\rho}: \pi_1(F,q) \rightarrow
PSL(2,\mc)$ in the conjugacy class of $\rho$. Indeed, we can push a
little $\gamma$ in the interior of the cylinder $\Cc$, and then use
the encoding of $\Rr(g,r)$ via $\Ii$-triangulations of $\mathcal{C}$
to check the claim. The point is that the very definition (\ref{ideq2defsec})
of the maps $W_a^I$ implies that the parallel transport operator along
$\gamma$ is not altered when we perturb it this way.


\section{The QHFT bordism category, and $\Ii$-cusps}\label{BORDCAT}
In this section we define the (2+1)-bordism category at the
basis of the QHFTs. We start with a naked, ``coordinate free''
category, then we progressively reach the final elaborated
marking. This shall incorporate the phase space parameters discussed
in Sections \ref{Dpar}-\ref{Ipar}. Finally, we introduce the notion of
$\Ii$-cusps.

\subsection{The (2+1)-bordism category}\label{bord+}

\paragraph{Naked bordism category.}\label{nakedcat}
Like in Section \ref{PARAM}, for every $(g,r)\in \N \times \ N$, such
that $g\geq 0$, $r>0$, and $r>2$ if $g=0$, we fix a compact closed
oriented {\it base} surface $S=S_g$, of genus $g$ with a set of $r$
marked points $V=V_{g,r}$. We denote by $-S$ the same surface with the
opposite orientation.  An {\it elementary object} of our {\it naked
category} either is the empty set, or is represented by a
diffeomorphism $\phi: \pm S \to \Sigma$. In other words it is a
parametrized surface $\Sigma$. Both the orientation and the marked
points $V$ transit from $\pm S$ to $\Sigma$ via $\phi$. Later the
surface $\pm S$ shall be equipped with further extra-structures such
as triangulations; we always stipulate that these extra-structures transit from
from $\pm S$ to $\Sigma$ via $\phi$. 

The pairs $(\pm S,\phi)$ are
considered up to the following equivalence relation: $(\pm S,\phi_1)$
is identified with $(\pm S,\phi_2)$ (i.e. they represent the same
elementary object) iff there is an oriented diffeomorphism $h:
\Sigma_1 \to \Sigma_2$, such that $(\phi_2)^{-1}\circ h \circ \phi_1$
pointwise fixes $V$ and is isotopic to the identity automorphism of
$S$ relatively to $V$. An {\it object} is a finite union of elementary
ones, where $(g,r)$ varies.

We define now the bordisms between objects, i.e. the {\it morphisms}
of the naked category.  Let $Y$ be an oriented compact $3$-manifold
with (possibly empty) boundary $\partial Y$. It is given a {\it input
  vs output} bipartition of the boundary components so that $\partial
Y = \partial_-Y \cup \partial_+ Y$. We can imagine that $\partial_- Y$
is ``at the bottom'' of $Y$, while $\partial_+Y$ is ``on the top''.
Each boundary component inherits the boundary orientation, via the
usual convention ``last is the ingoing normal''. We assume also that
it is given a properly embedded {\it link} $L\subset Y$, considered up
to proper ambient isotopy. Sometimes it is convenient to look at $L$
as $L= L_i\cup L_b$, where $L_i$ is the {\it internal} part of $L$
made by its closed connected components, while $L_b$ is the union of
the components homeomorphic to the interval $[0,1]$, with end-points
at some boundary components of $Y$ (possibly the same), and transverse
to $\partial Y$. For every boundary component $\Sigma$ of $Y$, we
assume that $|L_b\cap \Sigma|>0$, and that $|L_b\cap \Sigma|>2$ if
$g(\Sigma)=0$. Note that we do {\it not} require that every component
of $L_b$ connects $\partial_-$ with $\partial_+$. As the base surfaces
$S_g$ and the boundary components $\Sigma$ of $Y$ have given
orientations, we need to specify the ``sign'' of an object $[\phi: \pm S
\to \Sigma]$. 

\noindent Hence we can associate objects $\alpha_\pm$ to both
$\partial_\pm$.  We get in this way the bordism from $\alpha_-$ to
$\alpha_+$ {\it with support $(Y,L)$}.  We also allow that $\partial Y
= \emptyset$, i.e. $Y=W$ is a closed manifold, and $(W,L)$ is a
morphism from the empty set to itself. We stress that $L$ is non empty
in any case.

\paragraph {Introducing framings.}\label{frame}
For every $(S,V)$ as above, we introduce a {\it framing} at each
marked point $v\in V$. This means that we fix a system of disjoint
embedded segments $a_v$ in $S$ having the $v$'s as ``first''
end-point. We denote by $v''$ the other end-point. We adapt the above
definition of the {\it objects}, by requiring that $(\phi_2)^{-1}\circ
h \circ \phi_1$ is the identity on the $a_v$'s, and the isotopies are
relative to them.

We assume now that the above link $L\subset Y$ is {\it framed}, and we
denote it by $L_\Ff$. This means that $L_\Ff$ is a disjoint union of
properly embedded orientable {\it ribbons}. Each interior component of
$L_\Ff$ is homeomorphic to the annulus $S^1 \times [0,1]$, the other
components are homeomorphic to $I\times [0,1]$. On the boundary of
each ribbon we keep track of a {\it core} line of the form $X\times
\{0\}$ (here $X=S^1$ or $I$ resp.), for the corresponding component of
the unframed link $L$, and there is a {\it longitudinal} line $X\times
\{1\}$ that specifies the framing of the normal bundle of the parallel
core line.  This induces on each boundary component $\Sigma$ of $Y$ a
system of framed marked points (i.e.  $L\cap \Sigma$ is framed by
$L_\Ff \cap \Sigma$). Bordisms supported
by $(Y,L_\Ff)$ are defined similarly as above.

It is convenient to reformulate the bordism category with framed
links in an equivalent but slightly different way, which is closer to the
phase space parameters set up.

\paragraph {Zipping-unzipping.}\label{zip}
Fix a mid-point $v'$ on each arc $a_v$.  Then, let us {\it unzip} (cut
open) each $a_v$ at the open sub-interval $(v,v')$.  In this way we
get from $\pm S = \pm S_g$ an oriented surface $\pm F = \pm F_{g,r}$
with $r$ boundary components. Each boundary component of $F$ is a
bigon with vertices $v$ and $v'$; $v'$ is connected to $v''$ by the
sub-interval $[v',v'']$ of $a_v$. We can use these $\pm F$'s as
sources of elementary objects $[\phi: \pm F \to \Sigma]$.  Naturally,
also $\Sigma$ has now $r$ boundary components. To define equivalent
parametrized surfaces, we use (homotopy classes of) the oriented
diffeomorphisms of $\pm F$ that are the identity on the boundary.

Consider now $(Y,L_\Ff)$ as above. For every ribbon component of
$L_\Ff$ take a {\it mid-line} corresponding to $X \times \{1/2\}$.
This eventually gives us a triple
$\bar{\Lambda}=(\lambda,\lambda',\lambda'')$ of parallel {\it
  unframed} links in $Y$: the core line $\lambda=X\times \{0\}$, this
just introduced mid-line $\lambda'=X \times \{1/2\}$, and the
longitudinal boundary line $\lambda''=X \times \{1\}$. The trace of
$(L_\Ff,\bar{\Lambda})\cap \partial Y$ at each component $\Sigma$ of
$\partial Y$ makes a system $\{ a_u \}$ of disjoint segments; each one
has a marked end-point $u$, a mid-point $u'$ and the other end-point
$u''$.  Let us unzip each ribbon band at the open sub-band $X \times
(0,1/2)$. We get in this way a $3$-manifold ``with corners''
$\widetilde{Y}$. Its boundary $\partial \widetilde{Y}$ has two {\it
  horizontal} parts $\partial_\pm \widetilde{Y}$ contained in
$\partial_\pm Y$, and a {\it tunnel} part $\widetilde {L}_\Ff$. The
horizontal parts intersect the tunnel part at the corner locus; this
is a union of bigons contained in $\partial Y$. Each boundary
component of $Y$ corresponds to a horizontal boundary component of
$\widetilde{Y}$, still denoted by $\Sigma\in \partial_\pm$. Each {\it
  internal} tunnel component is homeomorphic to the torus
$\T=S^1\times S^1$; the other tunnel components are homeomorphic to
$\A=S^1 \times I$. Every tunnel component is made by the union of two
copies of $X \times (0,1/2)$, glued each to the other at $\lambda \cup
\lambda'$.

The horizontal boundary components can be considered as targets of
elementary objects $[\phi: \pm F \to \Sigma]$, each
$(\widetilde{Y},\widetilde{L}_\Ff)$ can be considered as the support
of a morphism between such objects.  All this is straightforward.
Clearly, we can {\it zip} back $(\widetilde{Y},\widetilde{L}_\Ff)$ to reobtain the initial $(Y,L_\Ff)$, so we have two equivalent
settings to describe the same bordism category.

\paragraph{The QHFT bordism category.}\label{dqftbordcat}
To define the (elementary) objects of our final bordism category, for
every $F$ as above, we fix as part of a marking an $e$-triangulation
$(T,b)$ of $F$, and a cocycle $z\in Z_I(T,b,x,y_x)_\beta$ (for some
$x$). So an elementary object is of the form $[(\pm
F,(T,b),x,z),\phi]$; we keep $(T,b)$, $x$ and $z$ fixed when we define
equivalent marked surfaces. Let $(\widetilde{Y},\widetilde{L}_\Ff) $
be a naked bordism as above, and assume furthemore that the boundary
objects $\alpha_\pm$ are equipped to be objects of the present QHFT
bordism category. Assume also that it is given a conjugacy class
$\rho$ of $PSL(2,\C)$-valued representations of $\pi_1(\widetilde{Y}
\setminus \widetilde{L}_\Ff)$. Then $(\widetilde{Y},\widetilde{L}_\Ff,
\rho)$ is the support of a QHFT bordism from $\alpha_-$ to $\alpha_+$
iff, for every elementary boundary object $[(\pm F,(T,b),x,z),\phi]$
we have $\phi^*(\rho)= [z]$.

\paragraph {Bordism Composition.}\label{bordcomp}
Consider a QHFT bordism $\Bb$ from $\alpha_-$ to $\alpha_+$, with
support $(\widetilde{Y},\widetilde{L}_\Ff,\rho)$, and another bordism
$\Bb'$ from $\alpha'_-$ to $\alpha'_+$, with support
$(\widetilde{Y}',\widetilde{L}'_\Ff,\rho')$. Assume that $\beta_+$ and
$\beta'_-$ are sub-objects of $\alpha_+$ and $\alpha'_-$ respectively,
such that $\beta_+ = - \beta'_-$. Then we can {\it glue} the two
bordisms at the common sub-objects. We get a new bordism (morphism) $\Bb'' :=
\Bb'*\Bb$ with support $(\widetilde{Y}'',\widetilde{L}''_\Ff,\rho'')$, from $\alpha''_-$ to
$\alpha''_+$, where $\alpha''_- = \alpha_- \cup (\alpha'_- \setminus
\beta'_-)$ and $\alpha''_+ = \alpha'_+ \cup (\alpha_+ \setminus
\beta_+)$. We say that $\Bb''$ is the {\it composition} of the bordism $\Bb$
followed by the bordism $\Bb'$.

\subsection {Cusps}\label{CUSP}
The cusps that we are going to introduce will allow us, in Section
\ref{TENSOR}, to turn the tunnel part of QHFT bordisms supported by a
triple $(\widetilde{Y},\widetilde{L}_\Ff,\rho)$ into ``toric ends'' or
``annular ends'' at infinity, thus replacing the corresponding
markings with the weaker, purely geometric, dependence w.r.t. the
$\rho$-holonomy at the meridians of $L_\Ff$.

\paragraph{Triangulated cusps.}\label{Tcusp}
We consider topological oriented cusps ${\bf C} = {\bf C}_\A = \A
\times [0,+\infty[,\  {\bf C}_\T = \T \times [0,+\infty[$ 
with {\it annular base} $\A=S^1\times
[-1,1]$ and {\it toric base} $\T=S^1\times S^1$ respectively. 

We fix first a specific class of branched triangulations of the bases
$\A, \T$. The idea is that of taking a branched triangulation
$(T',b')$ of the band $B=[0,1]\times X$, with $X= [-1,1]$ and $X= S^1$
respectively, having the vertices on the line $L=\{0,1\}\times X$;
then we {\it unzip}, i.e.  cut open, $B$ along $(0,1)\times X$, thus getting
a triangulation $(T,b)$ of $\A$ or $\T$ made by two copies of
$(T',b')$ that coincide on the line $L$. More precisely $B$ is
subdivided by a certain number of quadrilaterals $R$ having parallel
``vertical'' sides on $L$ and two interior parallel ``horizontal''
sides. Parallel sides of each $R$ have parallel orientation. Finally
a branched triangulation $(T',b')$ of $B$ is obtained by introducing
(in some way) an oriented diagonal on each $R$. In Fig. \ref{carre}
we see a fundamental domain of the base with a triangulation $(T,b)$ obtained by
starting with $2$ quadrilaterals.

\begin{figure}[ht]
\begin{center}
\includegraphics[width=4cm]{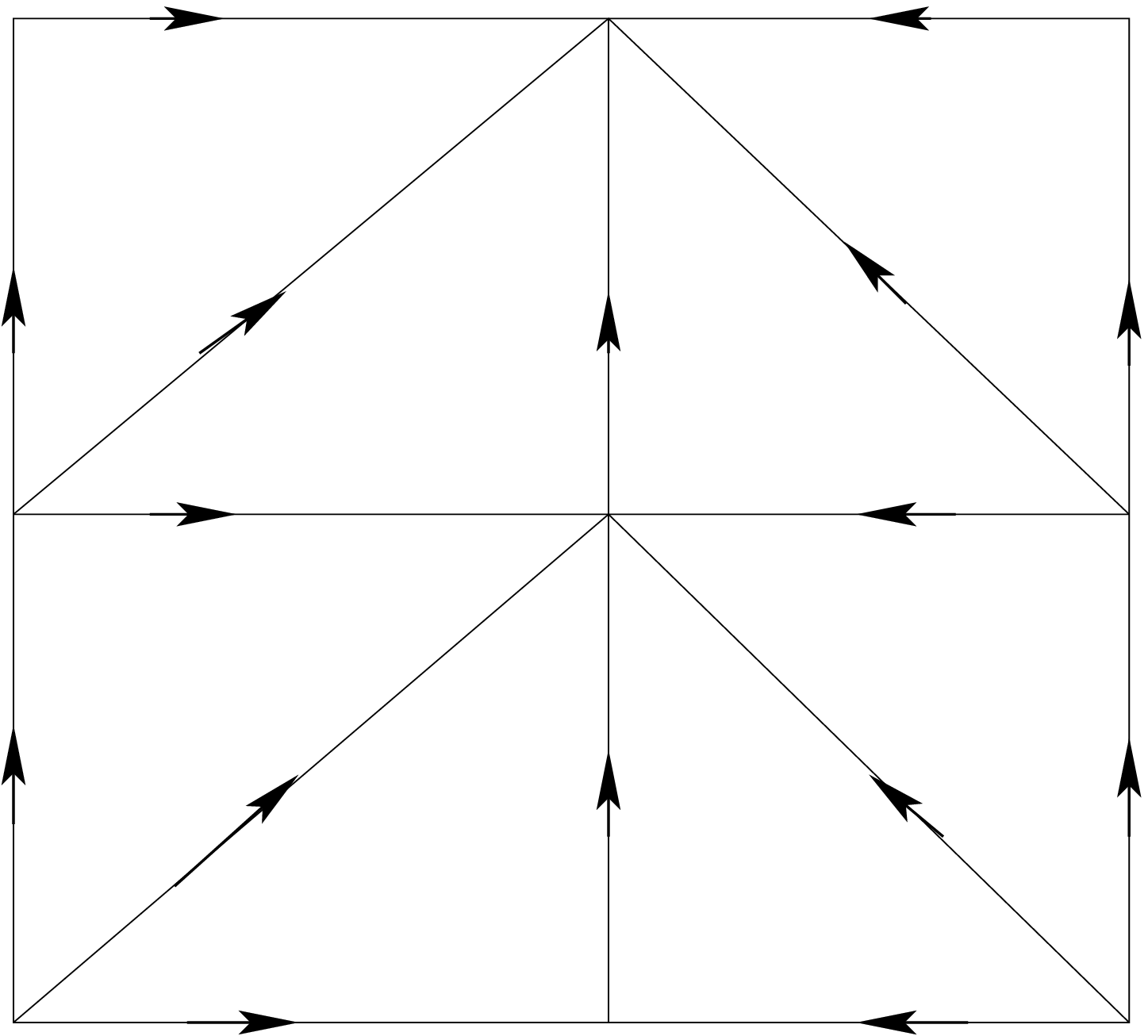}
\caption{\label{carre} A special triangulation of the torus.} 
\end{center}
\end{figure}

We consider now the one point-compactification {\bf C}$^*$
of a given cusp {\bf C}. From $(T,b)$ we get 
the a branched triangulation $(T^*,b^*)$ 
of {\bf C}$^*$, just by taking the one-point compactification of 
$(T,b)\times [0,\infty[$; more precisely we stipulate that the point
$\infty$ is a common vertex of all tetrahedra of  $(T^*,b^*)$, and 
is the opposite vertex to each triangle of $(T,b)$. The branching $b^*$
extends $b$, by imposing that $\infty$ is a pit for every 
branched tetrahedron. In Fig. \ref{tcusp} we see an example of
cusp triangulation (starting with one quadrilateral). 
  
\begin{figure}[ht]
\begin{center}
\includegraphics[width=5cm]{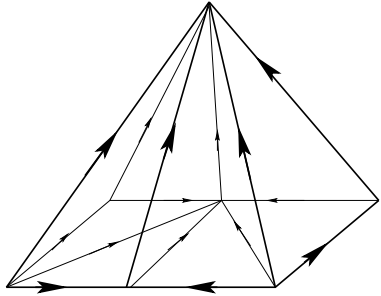}
\caption{\label{tcusp} An example of triangulated cusp.} 
\end{center}
\end{figure}

\paragraph{$\Ii$-cusps.}\label{Icusp}
Let $z$ be any {\it idealizable} $PB(2,\C)$-valued cocycle (i.e. with
values in the Borel subgroup of $PSL(2,\C)$) on some $(T,b)$ as above.
The idealization of each triangle is an ideal triangle with $3$
distinct ordered vertices belonging to $\C \subset \C\cup \{\infty
\}=\mc\mathbb{P}^1$. We look at such an ideal triangle as the base of
an ideal hyperbolic tetrahedron which as the further vertex at
$\infty$, and is oriented like the corresponding tetrahedron of
$(T^*,b^*)$. Doing it for every triangle, we obtain, by definition,
the {\it idealization} of $[T^*, b^*, z]$, which is called an
$\Ii$-{\it cusp}. We have in fact an $\Ii$-triangulation, that is:

\begin{lem} The cross-ratio moduli of every $\Ii$-cusp
verify the edge compatibility condition (\ref{ideq1defsec}) at each interior
``vertical edge'' (i.e. not contained in a base triangle).
\end{lem}

\noindent {\it Proof.} The cross-ratio moduli define a similarity
structure for each triangle of $T$. These extend to a {\it
  pseudo-similarity} structure of the whole base (see above the
  discussion of
the idealization procedure), with holonomy given
by the values of $z$ on simplicial representatives of generators of
  the fundamental group, or as at the end of Section \ref{Ipar}. In particular,
the holonomy is trivial on small closed loops on $T$ winding once
about the endpoint of an interior vertical edge of the
  $\Ii$-cusp. Hence (\ref{ideq1defsec}) is satisfied. \hfill$\Box$

\medskip


Finally we can precise the vertical relation over any $v$-vertex used
in Subsection \ref{Ipar}, before Proposition \ref{Ipar1}.
Given $z\in \Ii_a(T,b,x,y_x)_\beta$, the vertical tube over the
boundary loop associated to $v$ can be considered as
the base of an annular $\Ii$-cusp. Glue this $\Ii$-cusp to
the tube, so that the vertical edge over $v$ becomes interior. Then
the corresponding vertical relation is just an usual edge compatibility
condition.


\section {QHFT morphisms}\label{TENSOR} Consider a morphism $\Bb$ of the
QHFT bordism category from the object $\alpha_-$ to the object
$\alpha_+$, with support $(\widetilde{Y},\widetilde{L}_\Ff,\rho)$, as
in the previous Section.  For every odd integer $N\geq 1$, we will
show how to associate explicitely to $\alpha_\pm$ a finite dimensional
complex linear space $E_N(\alpha_\pm)$, and to the bordism a linear
morphism
$$
\Hh_N(\Bb) : E_N(\alpha_-) \to E_N(\alpha_+)$$
in a functorial way
w.r.t. the gluing of bordisms and morphism composition. This is
the technical core of the construction of the QHFTs. As usual we
associate to the empty object the ground field $\C$.

\subsection{Matrix dilogarithms and trace tensors}\label{MATTENS}

\paragraph{Review of the matrix dilogarithms.}\label{matdilog} The
building blocks of the QHFT morphisms are the {\it matrix
  dilogarithms} that we have introduced and studied in \cite{BB2}.
Here we limit ourselves to recall some of their qualitative
properties, what is just enough to follow the logic of the
construction.es
es

\smallskip

An $\Ii$-{\it tetrahedron} $(\Delta,b,w)$ consists of an {\it
oriented} tetrahedron $\Delta$ equipped with a {\it branching} $b$,
and a {\it modular triple} $w=(w_0,w_1,w_2)=(w(e_0),w(e_1),w(e_2))
\in (\C \setminus \{0,\ 1 \})^3$ such that (indices mod($\mz/3\mz$)):
$$w_{j+1} = 1/(1-w_j).$$
\noindent Hence $w_0w_1w_2 = -1$, and this gives a \emph{cross-ratio
  modulus} $w(e)$ to each edge $e$ of $\Delta$ by imposing that
$w(e)=w(e')$. We have already used such notions in Subsection
\ref{Dpar}. We will use the notations and conventions established
there.

Given any $\Ii$-tetrahedron $(\Delta,b,w)$, we consider an
extra-decoration made by two $\Z$-valued functions defined on the
edges of $\Delta$, called {\it flattening} and {\it integral charge}
respectively. These functions share the property that {\it opposite edges
take the same value}, hence it is enough to specify their values on the
edges $e_0,e_1,e_2$. We denote by $\log$ the standard branch of
the logarithm which has the arguments in $]-\pi,\pi]$.
\smallskip

For every $f=(f_0,f_1,f_2)$ with $f_i=f(e_i)\in \Z$, set
\begin{equation}\label{logb1}{\rm l}_j = {\rm l}_j(b,w,f)=\log(w_j) +
  i\pi f_j
\end{equation}
for $j=1$, $2$, $3$. We call ${\rm l}_j$ a \emph{log-branch} of
$(\Delta,b,w)$ for the edge $e_j$. We say that $(f_0,f_1,f_2)$ is a
\emph{flattening} of $(\Delta,b,w)$ if
$${\rm l}_0 + {\rm l}_1 +{\rm l}_2 = 0.$$ 

An {\it integral charge} of $(\Delta,b)$ is a function
$c=(c_0,c_1,c_2)$ with $c_i=c(e_i)\in \Z$, such that $c_0+c_1+c_2=1$. An $\Ii$-tetrahedron
endowed with a flattening and an integral charge is said {\it
flat/charged}. 
\medskip

For every $N>0$, any map
$$A: \C\setminus \{0,1\}\to {\rm Aut}(\C^N\otimes \C^N)$$
can be interpreted as a function of $\Ii$-tetrahedra via the formula:
$$A(\Delta,b,w):= A(w_0)^{*_b}.$$
Namely, put the standard tensor product basis on $\C^N\otimes \C^N$, so that $A=A(x)\in {\rm
Aut}(\C^N\otimes\C^N)$ is given by its matrix elements
$A^{\delta,\gamma}_{\beta,\alpha}$, where $\alpha, \ldots, \delta \in
\{0,\ldots,N-1\}$. We denote by $\bar{A}=\bar{A}(x)$ the inverse of
$A(x)$, with entries $\bar{A}^{\beta,\alpha}_{\delta,\gamma}$. The branching $b$ selects $w_0$ among the triple of
cross-ratio moduli. We use it also to associate to each
$2$-face of $\Delta$ one index among $\gamma,\delta,\alpha,\beta$. The
rule is shown in Fig. \ref{2CQDidealtensor}.

\smallskip

The matrix dilogarithm of rank $N$, $N\geq 1$ being any odd integer,
is an explicitely given ${\rm Aut}(\C^N\otimes \C^N)$-valued function
\begin{equation}\label{defRN}
\Rr_N(\Delta,b,w,f,c)= \Rr_N(w_0,f,c)^{*_b}
\end{equation}
defined on flat/charged $\Ii$-tetrahedra. The explicit formula is
given in (\ref{symRmatdil})-(\ref{symqmatdil}). Note that the matrix elements are
holomorphic w.r.t. the log-branches (up to sign when $N=1$); in particular, they depend on the
whole decoration, not only on $w_0$, so flattenings and charges are
incorporated in the above identification between tensors and decorated
tetrahedra. 

Each matrix dilogarithm $\Rr_N$ satisfies a finite set of fundamental
{\it five term identities}. These identities are supported by suitable
$\Ii$-flat/charged versions, called {\it transit configurations}, of
the basic $2\to 3$ bistellar move on $3$-dimensional triangulations
(sometimes called Pachner or Matveev-Piergallini move). This bare move
is shown on the top row of Fig. \ref{CQDfigmove1}. We postulate that
all the $5$ tetrahedra involved in the move are oriented and that they
induce opposite orientations on every common $2$-face. Hence, we have
two triangulations $T$ and $T'$ (by $2$ and $3$ tetrahedra resp.) of a
same oriented polyhedron, and each tetrahedron inherits the induced
orientation. The $2\to 3$ transit configurations involve an
appropriate procedure to transfer branchings, moduli, flattenings and
integral charges from the tetrahedra of $T$ to those of $T'$.

\begin{figure}[ht]
\begin{center}
\includegraphics[width=8cm]{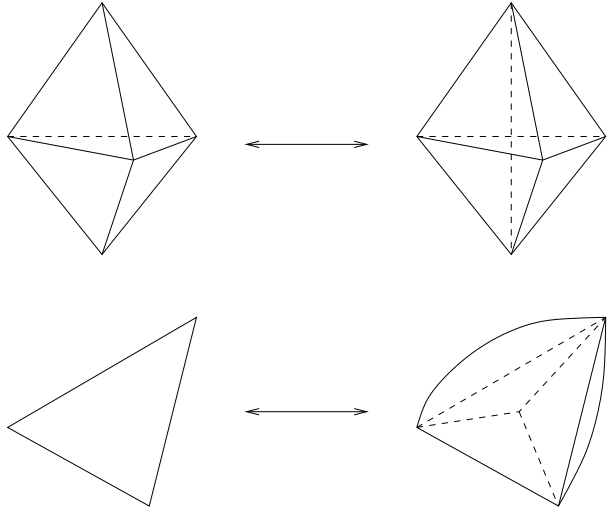}
\caption{\label{CQDfigmove1} The moves between singular triangulations.} 
\end{center}
\end{figure}
 
 The matrix dilogarithms satisfy also other relations corresponding to transit
 configurations associated to few other local $3$-dimensional
 triangulation moves, such as the so called {\it bubble move} which
 increases by one the number of interior vertices (see the bottom row
 of Fig. \ref{CQDfigmove1}).
 
 The matrix dilogarithms, as well as the elaborated extra decoration
 on $\Ii$-tetrahedra, arise from the solution of a {\it symmetrization
   problem} for a family of {\it basic} matrix dilogarithms, which
 satisfy only {\it one} peculiar five term identity supported by the
 top row move of Fig. \ref{CQDfigmove1}, called the {\it matrix
   Schaeffer's identity}.  This identity is characterized by
 determined geometric constraints on the cross ratio moduli. In the
 ``classical'' case $N=1$, the basic dilogarithm coincides with (the
 exponential of) the classical {\it Rogers dilogarithm}. The quantum ($N>1$) basic dilogarithms are
 derived from the $6j$-symbols for the cyclic representation theory of
 a Borel quantum subalgebra $\Bb_\zeta$ of $U_\zeta(sl(2,\mc))$, where
 $\zeta=\exp(2i\pi/N)$ (see \cite{K2}, Section 8 of \cite{BB2} or the
 Appendix of \cite{BB1}).

\paragraph{Trace tensors.}\label{trace}
Assume that we dispose of an $\Ii$-triangulation $(T,b,w)$ for a
compact oriented $3$-manifold $M$ (see Section \ref{Ipar}). Remind
that this means that the edge compatibility condition
(\ref{ideq1defsec}) is satisfied at each interior edge. Assume also
that every (abstract) $\Ii$-tetrahedron $(\Delta, b, w)$ of the
triangulation has a flat/charge $(f,c)$. Then, for every fixed odd integer
$N$, we can associate to each such $(\Delta, b, w, f,c)$ the
corresponding matrix dilogarithm $\Rr_N(\Delta, b, w, f,c)$.

\noindent A {\it $N$-state} of
$(T,b,w,f,c)$ is a function which associate to every triangle of the
$2$-skeleton of $T$ a value in $\{0,\dots,N-1 \}$. So, every $N$-state
determines indeed a matrix element of each matrix dilogarithm. As two
tetrahedra induce opposite orientations on any common face, our
identification rules (\ref{defRN}) together with Fig.
\ref{2CQDidealtensor} imply that an index at a common face is ``down''
for the $\Rr_N$ of one tetrahedron while it is ``up'' for the other.
By applying Einstein's rule of ``summing on repeated indices'' to the
matrix elements selected by each possible $N$-state of $T$, we get the
contraction (i.e. the {\it trace}) of these pattern of tensors
$\Rr_N(\Delta, b, w, f,c)$. We denote this trace by
\begin{equation}\label{tracedef}
\prod_{\Delta \subset T} \Rr_N(\Delta,b,w,f,c).
\end{equation} The type of
the trace tensor depends on the $b$-sign $\sigma$ of the boundary
triangles of $(T,b)$. The construction of $\Hh_N(\Bb)$ shall result
from a specific implementation of this procedure of taking the trace,
that includes suitable {\it global} constraints on flattenings and
integral charges.

\subsection{Bordism globally flat/charged $\Ii$-triangulations.}\label{globalfc}  
Let $\Bb$, $\alpha_\pm$, and $(\widetilde{Y},\widetilde{L}_\Ff,\rho)$ be
as at the beginning of this Section. Recall from Section \ref{BORDCAT}
that we have associated to $\Bb$ the pairs $(Y,L_\Ff)$,
$(Y,\bar{\Lambda})$ and $(Y,L)$. Here $L\subset Y$ is the initial
unframed link, the ribbon link $L_\Ff$ encodes a framed version of
$L$, $\bar{\Lambda} = (\lambda,\lambda',\lambda'')$ is made by $3$
unframed parallel copies of $L$, which lie on
$L_\Ff$. The links $\Lambda$ and $\Lambda'$ are on the tunnel boundary of
$\widetilde{Y}$, while $\Lambda''$ is transversal to the horizontal
boundary of $\widetilde{Y}$ at the image of the $v''$-vertices.

The first step consists in taking a so-called distinguished $\Dd$-triangulation
for $\Bb$, that we are going to define. For every elementary object
$[(F,(T,b),x, z), \phi]$ of $\alpha_-$ we take some idealizable
extension $\Cc^-(z,a)$ of $z$ to $\Cc(T,b)$ (see Section
\ref{Ipar}). Then we extend $\phi$ to a parametrization $\Phi$ of a
collar in $\widetilde{Y}$ of the corresponding horizontal boundary
component $\Sigma =\phi(F)\subset \partial_- Y$. For the elementary
objects in $\alpha_+$ of the form $[(-F,(T,b),x,z),\phi]$, we do
similarly, but we use $\Cc^+(z,a)$, that is the idealizable extensions
of $z$ obtained by exchanging the role of $F_-$ and $F_+$. Hence we
work with $0$-cochains supported at $\partial_-
\Cc(T,b)$ instead of $\partial_+ \Cc(T,b)$, as we did to define
$\Cc^-(z,a)$. By doing this for all elementary objects, we get a {\it
collar} $\Dd$-triangulation, that is a branched triangulation of a
collar of the whole horizontal boundary of $\widetilde{Y}$, equipped
with a $1$-cocycle $z_{coll}$ that represents the restriction of
$\rho$ to this collar.
\begin{defi}\label{Dtriang}{\rm
A {\it distinguished $\Dd$-triangulation} for $\Bb$ consists of a
$4$-uplet $\Kk =
(K,\bar{H},b,z)$ where:
\smallskip

(a) $(K,b)$ is a branched triangulation of $\widetilde{Y}$ that extends
a collar $\Dd$-triangulation.

\smallskip
(b) $(K,b)$ induces on each tunnel boundary component of
$\widetilde{Y}$ a triangulation of the type specified in Section
\ref{CUSP} for the cusp bases.
\smallskip

(c) The link $\bar{\Lambda}=\lambda \cup \lambda'\cup \lambda''$ is
realized by the sub-complex $\bar{H}= H\cup H'\cup H''$ of $K$, and
$\bar{H}$ is {\it Hamiltonian}, that is it contains all the vertices
of $K$. Because of (a) and (b), this is equivalent to require that
$H''$, which lies in the interior of $K$, contains all the {\it internal} vertices of $K$.  \smallskip

(d) $z$ is a $PSL(2,\C)$-valued $1$-cocycle on $(K,b)$ such that: $z$ extends the collar cocycle $z_{coll}$; $\rho = [z]$ on the whole $\widetilde{Y}$; $z$ is idealizable; the restriction of $z$ to each tunnel component of $\partial
\widetilde{Y}$ is of the kind specified in Section \ref{CUSP} for the
bases of $\Ii$-cusps.}
\end{defi}

The next step consists in passing from the above $\Dd$-triangulation
$\Kk=(K,\bar{H},b,z)$ to an $\Ii$-triangulation
$\Kk_\Ii=(K_I,\bar{H},b,w)$ for $\Bb$.  To do it, Consider every tunnel
boundary component as the base of a topological cusp, oriented in such
a way that, by glueing these cusps, we get the oriented $3$-manifold
$Y\setminus L$. Now take the  idealization of $z$, and glue the corresponding $\Ii$-cusps (see Section \ref{CUSP}). 
The resulting family of $\Ii$-tetrahedra defines an
$\Ii$-{triangulation} $\Kk_\Ii$ for $\Bb$; at every internal edge it satifies the
edge compatibility condition (\ref{ideq1defsec}).

We specify now which kind of {\it global} flattening and integral
charge are carried by $\Kk_\Ii=(K_I,\bar{H},b,w)$. For every
elementary object $[(\pm S,(T,b),x,z),\phi]$ in $\alpha_\pm$, and for
every edge $e$ of $(T_\pm,b_\pm)$, we have the
values $W_\pm(z)(e)$ of the $\Ii$-boundary function defined in
Subsection \ref{Dpar}.  We can extend the definition of $W_\pm(z)$ to
the edges in the tunnel boundary components of $K_I$. By taking the union over the
elementary objects we get the {\it $\Ii$-boundary data} $W_\pm$ of the
objects $\alpha_\pm$.

Recall that we denote by log the standard branch of the logarithm with
arguments in $]-\pi, \pi]$.

\begin{defi}\label{flat/charge}{\rm 
A {\it global flattening} on  $(K_I,\bar{H},b,w)$ is a collection $f$
of flattenings for the $\Ii$-tetrahedra of $(K_I,\bar{H},b,w)$
such that:
\smallskip
 
(1) At each {\it internal} edge $e$ of $K_I$, the associated log-branches
formally satisfy the log of the edge compatibility condition
(\ref{ideq2defsec}), that is:

\begin{equation}\label{2ideq2defsec}
\sum_{h\in \epsilon^{-1}(e)} *\ {\rm l}(h)= 0 = \log (1)
\end{equation}
where $\epsilon$ is the identification map that
associates to every abstract edge its image in $K_I$ (via
the face pairings), and $*=\pm 1$ according to the
$b$-orientation of the tetrahedron that contains the abstract edge $h$.
\smallskip

(2) The condition (\ref{2ideq2defsec}) extends at every {\it boundary}
edge $e$ of $K_I$ by requiring that
\begin{equation}\label{Wideq2defsec}
\sum_{h\in \epsilon^{-1}(e)} *\ {\rm l}(h)= \log (W_\pm(e)).
\end{equation}
\medskip

A {\it global integral charge} on $(K_I,\bar{H},b,w)$ is a collection
of integral charges over the tetrahedra of $K_I$ such that:
\smallskip

(a) The sum of charges about every {\it internal} edge of 
$K_I \setminus \bar{H}$ is equal to $2$.
\smallskip

(b) The sum of charges at every {\it boundary} edge of 
$K_I$ is equal to $1$.

\smallskip

(c) The sum of charges about every edge of 
$H''$ equals $0$.
\smallskip

An $\Ii$-triangulation $(\Kk_\Ii,f,c)$ equipped with a global flattening
$f$ and a global integral charge $c$ is said {\it (globally) flat/charged}.  
}
\end{defi}

\paragraph{The cohomological charge.}  
Let $(\Kk_\Ii,f,c)$ be a flat/charged $\Ii$-triangulation for $\Bb$ as
above.  The reductions mod$(2)$ of both $f$ and $c$ induce cohomology
classes $[f],[c] \in H^1(\widetilde{Y};\Z/2\Z)$. Moreover, we have
integral classes $<f>,<c> \in H^1(\partial\widetilde{Y}; \Z)$ on
$\partial\widetilde{Y}$. These classes are defined as follows. For any
mod($2$) (resp. integral) $1$-homology class $a$ in $Y\setminus L$
(resp. $b$ in $\partial\widetilde{Y}$), realize it by a disjoint union
of (resp. oriented and essential) closed paths transverse to the
triangulation $\Kk_\Ii$ and `without back-tracking', i.e. such that it
never departs from a $2$-face of a tetrahedron (resp. $1$-face of a
triangle) by which it just entered. Then the mod($2$) sum of the
flattenings or charges associated to the angular sectors that we cross
when following such paths in the interior of $\widetilde{Y}$ define
the value of the corresponding class on $a$. Similarly, the {\it
  signed} sum of the log-branches or charges of the corners
that we cross when following such paths on $\partial
\widetilde{Y}$ define the value of the corresponding class on $b$; for
each vertex $v$ the sign is $*_b$ (resp. $-*_b$) if, with respect to
$v$, the path goes in the direction (resp. opposite to the one) given
by the orientation of $\partial \widetilde{Y}$.

\noindent The set $\{[f],[c],<f>,<c>\}$ is called the {\it
  cohomological charge} of $(\Kk_\Ii,f,c)$. It is a fact (see
\cite{BB2}) that this cohomological charge is preserved by the transit
configurations mentioned in Subsection \ref{MATTENS}. Let us denote by
$\partial_{tunnel}\widetilde{Y}$ the tunnel boundary
components of $\widetilde{Y}$. It is convenient to {\it normalize} the
cohomological charge by requiring that: 

\smallskip

(i) $<f>=<c>=k$, and $k$ is identically $0$ except possibly on $\partial_{tunnel}\widetilde{Y}$;

(ii) $[f]=[c]= h$, and $i^*(h) = k$ mod($2$), where $i: \partial_{tunnel}
\widetilde{Y} \rightarrow \widetilde{Y}$ is the inclusion map.

\smallskip

\noindent With this normalization the whole boundary configuration of a flat/charged
triple $(\Kk_\Ii,f,c)$ {\it is completely determined by the objects
  $\alpha_\pm$}, and the cohomological charge reduces to $(h,k)$. From
now on we consider only triples $(\Kk_\Ii,f,c)$ with such normalized
cohomological charge. The results below can be extended to the case
when the pairs $([f],<f>)$ and $([c],<c>)$ are independent, or not
equal to $0$ on $\partial\widetilde{Y} \setminus
\partial_{tunnel}\widetilde{Y}$; the only difference lies in the
constraints on the cohomological charge when gluing QHFT bordisms (see
Remark \ref{remclass}).

\subsection{The QHFT bordism tensors
  $\Hh_N(\Kk_\Ii,f,c)$}\label{triangtensor}

Fix an odd integer $N\geq 1$. Let $\Bb$, $\alpha_\pm$ and
$(\widetilde(Y),\widetilde{L}_\Ff,\rho)$ be as usual, and let
$(\Kk_\Ii,f,c)$ be a flat/charged $\Ii$-triangulation
with normalized cohomological charge for $\Bb$.  

First we specify the linear spaces $E_N(\alpha_\pm)$.  Recall that it
is defined a sign function $\sigma_\pm$ on the set of triangles of
$\alpha_\pm$.  Let us fix an ordering of the
elementary objects, together with an ordering of the triangles of each elementary
object. Thus we get a lexicographical order on the whole set of
triangles of $\alpha_\pm$. Fix a complex linear space $E=E_N$ of dimension
$N$, endowed with a given basis, so that it is identified with $\C^N$. Write
$E=E^1$ and $E^{-1}= E'$ for its dual space. Also $E'$ is canonically
identified with $\C^N$. Then set $E_N(\alpha_\pm)$ to be the tensor
product of the $|J_\pm|$ ordered spaces $E^{\sigma_\pm (t)}$, where
$t$ spans the set of triangles of $\alpha_\pm$ and $|J_\pm|$ is the
number of these triangles. The space $E_N(\alpha_\pm)$ is identified with the tensor
product of $|J_\pm|$ copies of $\C^N$.

We consider every matrix dilogarithm $\Rr_N \in {\rm Aut}(\C^N\otimes
\C^N)$ as an element of $(E_N')^{\otimes 2}\otimes (E_N)^{\otimes 2}$.
The trace tensor (\ref{tracedef}) for the pattern of
matrix dilogarithms associated to $(\Kk_\Ii,f,c)$ gives us a morphism
$$\Hh_N(\Kk_\Ii,f,c) \in {\rm Hom}(E_N(\alpha_-),E_N(\alpha_+)) \ .$$
We can state now the main technical result in the construction
of QHFT$_N$.
\begin{teo}\label{main} Let 
$\Bb$ be a QHFT bordism between objects $\alpha_\pm$. Then:
\smallskip

(1) Flat/charged $\Ii$-triangulations $(\Kk_\Ii,f,c)$ for $\Bb$, with
any prescribed normalized cohomological charge $(h,k)$, do exist.
\smallskip

(2) Let $(\Kk_\Ii,f,c)$ be such a triangulation with cohomological
charge $(h,k)$. Let $v_\partial$ and $v_I$ be respectively the number
of boundary and internal vertices of the triangulation.  Then, for
every odd integer $N\geq 1$, {\rm up to a sign and a $N$th-root of
  unity multiplicative factor}, the normalized linear map
$$N^{-(v_\partial/2+v_I)}\Hh_N(\Kk_\Ii,f,c) \in {\rm
  Hom}(E_N(\alpha_-),E_N(\alpha_+))$$ only depends on
the triple $(\alpha_\pm,h,k)$, so that the {\rm bordism tensor} 
$$\Hh_N(\Bb,h,k):=
N^{-(v_\partial/2+v_I)}\Hh_N(\Kk_\Ii,f,c)$$ is well defined (up to the
above phase ambiguity).
\end{teo}

\noindent {\it Proof.} The proof of this theorem is technically
demanding, but it follows strictly the arguments of \cite{BB1,BB2}. As
there are only slight differences, we limit ourselves to few comments.

\noindent Point (1). The existence of $\Dd$-triangulations $\Kk$ is
an essentially straightforward adaptation. We know that integral
charges and flattenings, with arbitrary cohomological charge, exist on
any closed $3$-dim. ($\Ii$-)triangulation. Moreover, they make an affine
space over a same integral lattice $\mathfrak{L}$, generated by
determined vectors attached to the abstract stars of all the edges
(see \cite[\S 4]{BB1} and \cite[\S 6]{BB2}). The lattice is fixed by
the cohomological charge. Consider the double $D\Kk_\Ii := \Kk_\Ii
\cup_\partial (-\Kk_\Ii)$. By symmetry of the triangulation about
$\partial \Kk_\Ii \subset D\Kk_\Ii$, the above facts imply that we can
find a flat/charge on $D\Kk_\Ii$ that induces one on $\Kk_\Ii$; in
particular (2) and (c) in Definition \ref{flat/charge} are satisfied.
An easy Mayer-Vietoris argument show that we recover the cohomological
charge $(h,k)$ by properly choosing it on $D\Kk_\Ii$ (this agrees with the normalization $<f>=<k>$). Finally, Definition
\ref{flat/charge} (2)-(c) just kill the generating vectors of
$\mathfrak{L}$ at the edges of $\partial \Kk_\Ii$. So the sets of
flattenings and integral charges on $\Kk_\Ii$ are affine spaces over a
lattice $\mathfrak{L}'$ generated by the interior edges of $\Kk_\Ii$
only.

\noindent For (2), the key point consists in showing that any two
arbitrary flat/charged $\Ii$-triangulations for $\Bb$ are connected
via a finite sequence of transit configurations (see Subsection
\ref{MATTENS}) in the interior of $\Kk_\Ii$, so that after each
transit we still have a flat/charged $\Ii$-triangulation of $\Bb$.
Such sequences preserve both the whole boundary configuration and the
normalized cohomological charge, and they allow also to vary the
flattenings and integral charges by any multiple of the generators of
the corresponding lattice $\mathfrak{L}'$.

\noindent Finally, the fundamental
identities satisfied by the matrix dilogarithms (the five term ones as
well as the ones supported by other triangulation local moves) imply
that the trace tensor is transit invariant, up to the specified phase
ambiguity. In particular, the normalization factor $N^{-v_I}$ is
needed for the bubble transit invariance (the factor
$N^{-v_\partial/2}$ is introduced to have a good behaviour with
respect to bordism composition, see Prop. \ref{functQHFT} below).

\noindent A
new ingredient w.r.t. to \cite{BB1, BB2} is the presence of
$\Ii$-cusps. The transits preserve the {\it germ at infinity} of each
cusp, that is the topological singularity corresponding to the
collapsed tunnel boundary components of $\widetilde{Y}$, and allow
also to modify the base triangulation. So the bordism tensors do not depend
on any fixed representative of $\rho$ at the link meridians.\hfill$\Box$ \medskip
 
\medskip

Recall the $\Ii$-boundary data $W_{\pm}$ before Definition \ref{flat/charge}. We have:

\begin{cor} The QHFT morphisms supported by a bordism between objects
  $\alpha_\pm$ define complex analytic functions (up to the phase
  ambiguity mentionned in Th. \ref{main}) on
  open dense subsets of the space of $\Ii\partial$-parameters $W_-
  \times W_+$ for $\alpha_+ \cup \alpha_-$.
\end{cor}

\noindent {\it Proof.} Given a QHFT bordism $\mathcal{B}$, fix an
arbitrary branched triangulation $(K,b)$ as in Definition
\ref{Dtriang}. The space of log-branches supported by $(K,b)$ is
linear, and the matrix dilogarithms (\ref{defRN}) depend
holomorphically on the log-branches (up to sign when $N=1$). By
definition, the QHFT tensors $\Hh_N(\Bb)$ are specific (normalized)
contractions of matrix dilogarithms, and Theorem \ref{main} (2)
implies that they are eventually functions of the right-hand side members of the
identities (\ref{Wideq2defsec}). As these are logarithms of
the $\Ii\partial$-parameters, we deduce that the QHFT tensors are indeed
analytic functions on a dense subset of $W_- \times W_+$.\hfill$\Box$

\medskip

To treat the morphism composition, we incorporate in the gluing set
up the fixed ordering on the set of triangles; moreover we include in
the definition of $Z_I(T,b,x,y_x)_\beta$ (see Section \ref{Dpar}-\ref{Ipar}) the
{\it genericity assumption} that for no boundary edge $e$ of
$(K_I,b)$ the $\Ii\partial$-parameter $W_\pm(e)$ belongs to the
negative real ray. This is a mild assumption, and it is completely
uninfluent if the surface triangulations are quasi-regular. The
advantage of this assumption is the following fact, which follows from
the proof of Lemma \ref{W}:
\medskip

\noindent If $\Bb'' = \Bb'*\Bb$ is a bordism obtained as in Section
\ref{BORDCAT}, and $(\Kk,f,c)$, $(\Kk',f',c')$ are flat/charged
$\Ii$-triangulations with normalized cohomological charges $(h,k)$ and
$(h',k')$, they automatically glue to give a flat/charged
$\Ii$-triangulation $(\Kk'',f'',c'')$ for $\Bb''$, for some $(h'',
k''):= (h'*h, k'*k)$.

\begin{remark} \label{remclass} {\rm Even if $h=h'=0$ it
    is possible that $h''\ne 0$. In general $h''$ depends on the
    whole pair $(\Bb,\Bb')$. This follows from an easy application of
    a Mayer-Vietoris argument. But if the glued part of the boundary
    is connected, or is a boundary in $\Bb * \Bb'$, then $h=h'=0$ implies
    $h''=0$.}
\end{remark}

The following functoriality result is a consequence of Theorem \ref{main}.

\begin{prop} \label{functQHFT} With the above hypothesis we have:
$$\Hh_N(\Bb'',h'',k'')=_N
\Hh_N(\Bb',h',k')\circ \Hh_N(\Bb,h,k)$$
where ``$=_N$'' denotes the phase ambiguity mentioned in Theorem
\ref{main} (2), and $\circ$ is the morphism composition.
\end{prop}

There is also a Hermitian property of QHFT morphisms (proved as for closed
manifolds, see Prop. 4.29 of \cite{BB2}):

\begin{prop}\label{unitarity} Write $\bar{\Bb}$ for the QHFT bordism
with opposite total space orientation for $\rho$ (i.e. opposite
orientation for $Y$ and complex conjugate holonomy $\bar{\rho}$). Then:
$$\Hh_N(\bar{\Bb},h,-k) =_N \Hh_N(\Bb,h,k)^{*T}$$
where $^{*T}$ denotes the matrix hermitian conjugation.
\end{prop}

\subsection{QHFT partition functions}\label{PARTF}

Here we discuss numerical invariants, called {\it QHFT partition
  functions}, for closed manifolds, i.e. manifolds with empty boundary.

\smallskip

Assume that $Y=W$ is closed, and that
it contains a {\it non empty} link $L_\Ff$. Then $(W,L_\Ff)$ supports
a QHFT bordism from the empty object to itself. The associated tensor
$\Hh_N(W,L_\Ff,\rho,h,k)$ is a scalar (well defined up to ``$=_N$'').
We can also set the cohomological charge to be $(0,0)$. So we get an
invariant $\Hh_N(W,L_\Ff,\rho)$ for the triple $(W,L_\Ff,\rho)$. Typical examples of such triples are given by cone hyperbolic manifolds
$W$ with cone locus an unframed link $L$. The character $\rho$ is just
the hyperbolic holonomy of $W \setminus L$, and we take an arbitrary framing for $L$.

Here is another way to express these partition functions, in terms of
manifolds with toric boundary. In that situation, fix an ordered basis
$(m_i,l_i)$ for the integral homology of each boundary torus. Then we
can construct the pair $(W,L_\Ff)$, where $W$ is obtained by Dehn
filling of $Y$ along the $m_i$'s, and $L_\Ff$ is the disjoint union of
the cores of the surgered solid tori, framed by the $l_i$'s. In this
way, the invariants $\Hh_N(W,L_\Ff,\rho)$ are viewed as invariants of
$(Y,\{(m_i,l_i)\}_i)$.

\smallskip

Let us introduce a QHFT bordism category that covers a {\it more
  restricted range} of topological/geometric situations. We propose
this variation because the corresponding partition functions exactly
equal the invariants for triples $(W,L,\rho)$ already constructed in
\cite{BB1, BB2}.

Consider the naked bordism category with {\it unframed} links
introduced at the beginning of Section \ref{BORDCAT}. We assume also
that $(Y,L)$ is equipped with a $PSL(2,\mc)$-character $\rho$ which is
defined on $\pi_1(Y)$; hence $\rho$ is trivial at the link meridians.
To define the elementary objects, we use triangulations $(T',b')$ of
$(S,V)$ of the type discussed in Section \ref{PARAM}. Moreover we
consider cocycles $z$ which have idealizable extensions $\Cc^-(z,a)$
to the (triangulated) cylinder $S\times [-1,1]$. This give us a notion
of ``QHFT$^0$'' bordism category. Note that the $\Ii$-cusps are no
longer present. We consider $\Dd$-triangulations for which the link
$L$ is realized as a Hamiltonian subcomplex, and we have only the mod
$(2)$ cohomological charge $h$. The analog of Theorem \ref{main}
associates tensors $H_N(\Bb,h)$ (still defined mod $(=_N)$) to every
QHFT$^0$ bordism $\Bb$. For triples $(W,L,\rho)$ with $W$ closed, we
obtain partition functions $H_N(W,L,\rho)$.

So, we dispose,for every $N$, of the two partitions functions
$\Hh_N(W,L_\Ff,\rho)$ and $H_N(W,L,\rho)$. The next result shows that
they are very similar indeed.

\begin{prop} We have $\Hh_N(W,L_\Ff,\rho)=_N
  H_N(W,\bar{\Lambda},\rho)$, where, as in Section \ref{BORDCAT}, the
  link $\bar{\Lambda}$ is made by $3$ specific parallel copies of $L$
  given by the framing $L_\Ff$.
\end{prop}

\noindent {\it Proof.} We consider supports that are very
close each to the other. Recall the {\it zipping} procedure from
the subsections \ref{e-triang} and \ref{bord+}. Since $\rho$ has trivial holonomy at the
meridians of $L_\Ff$, we can zip back any $\Dd$-triangulation of
$(W,L_\Ff,\rho)$, in
particular the cocycle, to
get a $\Dd$-triangulation used to compute
$H_N(W,\bar{\Lambda},\rho)$ (passing to its idealization).  After the
zipping there are no longer cusps, and the cross-ratio moduli of the
survived tetrahedra are unchanged.

\noindent So, to conclude we must check that the cusp contributions
to $\Hh_N(W,L_\Ff,\rho)$ are trivial. By looking at Figure \ref{carre}, we see
that each cusp tensor is supported by two patterns of flat/charged
$\Ii$-tetrahedra glued in a mirror-like fashion w.r.t. the central
vertical line. Mirror $\Ii$-tetrahedra have the same cross-ratio
moduli. So we may choose exactly the same flattenings (resp. integral
charges) for them; one example is where each right-angle in the figure
is given the flattening $\pm 1$ (according to the sign of minus the
imaginary part of the corresponding cross-ratio modulus), whereas the
acute angles are given the value $0$. For the charges we put
correspondingly the values $1$ or $0$. It is easily seen that such a
choice makes a flat/charged $\Ii$-triangulation, where the
cohomological charge is identically $0$. The branching induces
opposite orientations of the mirror patterns of flat/charged
$\Ii$-tetrahedra.  So the respective automorphisms are inverse one to
the other, and each cusp contribution to $\Hh_N(W,L_\Ff,\rho)$ is an identity tensor.\hfill$\Box$

\begin{remark}\label{linkre}{\rm It is known \cite {BB1, BB2} that in the
    ``classical'' case $N=1$, the link is immaterial in the
    construction of $H_N(W,L,\rho)$. But in the ``quantum cases''
    $N>1$, the presence of a non empty link is essential for the
    theory, and its choice alters the value of $H_N(W,L,\rho)$. Also,
    when $N=1$, it is known that $H_1(W,\rho)$ equals
    $\exp(1/i\pi({\rm CS}(\rho) +i {\rm Vol}(\rho)))$ (see
    \cite{BB2}, Section 6).}
\end{remark}

\subsection{QHFT supported by product bordisms}\label{MODREP}

In this section we consider the product bordisms $\widetilde{Y}=
F\times [-1,1]$, containing the tunnel boundary $\partial F \times
[-1,1]$, and equipped with a conjugacy class $\rho$ of
$PSL(2,\C)$-valued representations of $\pi_1(F)=\pi_1(\widetilde{Y})$.
This corresponds to triples $(Y,L_\Ff,\rho)$, where $Y=S\times [-1,1]$
and $L_\Ff$ is associated to $L= V\times [-1,1]$.

\noindent We are interested in the families of QHFT morphisms
supported by $\widetilde{Y}$, {\it when the marking of the boundary
  components $F \times \{-1\}$ and $F \times \{1\}$ vary}.  In
fact, we can consider separately the variations that
concern the parametrizing homeomorphisms, the $e$-triangulations, and
the cocycles, respectively. For each one, we can also restrict
to {\it elementary} generating variations such as Dehn twists, flips,
and gauge transformations induced by $0$-cochains supported at one
vertex, respectively.  \smallskip

Consider in particular objects of the form $[(\pm F,
(T,b),x,z),\psi]$, where we stipulate that $(\pm F, (T,b),x,z)$ {\it is
fixed and only the homeomorphism $\psi$ varies}.  Given $\psi_-$ and
$\psi_+$, set $\psi= \psi_+^{-1}\psi_-$ and $[\psi]$ for the
corresponding element in the mapping class group ${\rm Mod}(g,r)$.
\smallskip

Denote by $(\Tt_{[\psi]}, L_{[\psi]})$ the {\it mapping torus} of
$\psi$, which only depends on $[\psi]$.  Here $L_{[\psi]}$ is the
tunnel boundary of $\Tt_{[\psi]}$; its components are tori. Denote again $\rho$
the pull-back on $\Tt_{[\psi]}$ of the given conjugacy class $\rho$,
via the natural projection onto $F$. For the cylinders and the mapping
tori, the cohomological charge shall be trivial so we omit to mention
it. Fix an odd integer $N\geq 1$.  We denote by $q_N(\rho,[\psi])$ the
QHFT partition function for the triple $(\Tt_{[\psi]},
L_{[\psi]},\rho)$.  The cylinder $\widetilde{Y}$ only depends on the
pair $(g,r)$, and the objects only depend on $[\phi_\pm]$ (as we are
assuming that $(\pm F, (T,b),x,z)$ is fixed). So we write
$\Hh_N(r,g,[\psi_-],[\psi_+])$ for $\Hh_N(Y,L_\Ff,\rho)$. Set $d_{g,r}(N) = {\rm dim}
(E_N(\alpha_\pm))= N^{2(2g-2+3r)}$.

\begin{lem}\label{CYLtensor} 
(1) The QHFT tensor $\Hh_N(g,r,[\psi_-],[\psi_+])$ only depends on
$[\psi]$. So we denote it by $\Hh_N(g,r,[\psi])$.
\smallskip

(2) $\Hh_N(g,r,[id]) =_N {\rm Id}$, and ${\rm Trace}\ \Hh_N(g,r,[id])=_N d_{g,r}(N) =_N q_N(\rho,[id])$.
\smallskip

(3) $\Hh_N(g,r,[h_2])\circ \Hh_N(g,r,[h_1]) =_N
\Hh_N(g,r,[h_2h_1])$.  Hence, in particular:
$\Hh_N(g,r,[\psi^{-1}])\circ \Hh_N(g,r,[\psi]) =_N {\rm Id}.$
\smallskip

(4)  ${\rm Trace}\ \Hh_N(g,r,[\psi])=_N q_N(\rho,[\psi])$. 
\end{lem}      
\noindent{\it Proof.} By Theorem \ref{main} (2), the tensor $\Hh_N(g,r,[\psi_-],[\psi_+])$
is equal (up to the ``$=_N$'' ambiguity) to the one supported by the
mapping cylinder of $[\psi]$, where $F_-$ is parametrized by the
identity and $F_+$ is parametrized by $\psi$. Indeed, the
homeomorphism $\psi_-^{-1} \times {\rm id}$ relates the two supporting
cylinders. Point (1) follows. Also, this gives
$$\Hh_N^2(g,r,[id]) =_N \Hh_N(g,r,[id])$$
i.e. $\Hh_N(g,r,[id])$ is an
idempotent. By using a triangulation of the trivial mapping cylinder
of the form $\Cc^+(z,a)*\Cc^-(z,a)$ and the fact that the $\Rr_N$'s
are automorphisms, we see from (\ref{tracedef}) that $\Hh_N(g,r,[id])$
read as a tensor product of invertible matrices. Hence it is
invertible. The claims (3)-(4) are consequences of Prop.
\ref{functQHFT}. \hfill $\Box$

\medskip

If, with the usual notations, we vary now only the cocycles $z_\pm$ in
specifying the different objects, we still get {\it invertible}
tensors $\Hh_N(g,r,z_\pm)$ by the same arguments. If we change the
triangulation $(T_-,b_-)$, for instance, by a branched elementary
flip, getting $(T_+,b_+)$, there is a unique way to define $z_+$ in
such a way that it coincides with a given $z_-$ on the common edges.
Generically also $z_+$ belong to $Z_I(T_+,b_+,x,y_x)_\beta$, and again
we get an invertible representing tensor $\Hh_N(T_\pm,b_\pm,z_\pm)$.
The effect of such changes of the marking (for a given $\rho$) on a
bordism tensor $\Hh_N(g,r,[\psi])$ is a conjugation by $\Hh_N(g,r,z_\pm)$ or $\Hh_N(T_\pm,b_\pm,z_\pm)$. So we have:

\begin{cor}\label{mcgrep} 
  For every $N$, $\rho$, and $(T,b,z)$ as above, the tensors
  $\Hh_N(g,r,[\psi])$ realize a {\rm $d_{g,r}(N)$-dimensional linear
    representation} mod $(=_N)$ of the mapping class group ${\rm
    Mod}(g,r)$. Hence, by varying $(T,b,z)$, there is a well-defined
  conjugacy class of representations mod $(=_N)$ of ${\rm Mod}(g,r)$
  associated to QHFT$_N$.
\end{cor}

\subsection{Universal QHFT environment}\label{VARIATION}

Here we indicate the most general environment where the matrix dilogarithm
technology applies, and QHFT, {\it formally} at least, makes sense. 

\smallskip

We can associate trace tensors to any {\it roughly
flat/charged $\Ii$-triangulated} compact oriented {\it
pseudo-manifold} $Z$.  Such a triangulation is given, as usual, by a
family of oriented abstract flat/charged $\Ii$-tetrahedra equipped
with certain orientation reversing identifications of some $2$-faces,
in such a way that $Z$ is the quotient space. ``Pseudo-manifold''
means that the non-manifold points of $Z$ are contained in the set of
vertices of the triangulation. ``Roughly'' means that we require {\it neither}
that the edge compatibily condition is satisfied at the
internal edges, {\it nor} that the collection of local flattenings and
integral charges satisfy any global constraint. So we are in a
situation much more general than that of Section \ref{triangtensor}.

We consider these triangulated and decorated pseudo-manifolds {\it up
to the equivalence relation generated by the transit configurations
mentioned before Fig. \ref{CQDfigmove1}} (i.e. the flat/charged
$\Ii$-transits described in detail in \cite{BB1,BB2}).

The notion of boundary is well defined for this class of triangulated
pseudo-manifolds, so we can construct a consistent
bordism category. In fact, we can also relax the requirement that the
input/output bipartition of the boundary of the bordisms is made by
the disjoint union of {\it closed} boundary components. We can consider adequate
triangulated portions of the boundaries, such that the
category is stable under bordism composition (see for instance
\cite{T} for such a general notion of bordism category).

In this way, a single flat/charged $\Ii$-tetrahedron can be considered
as a bordism between the two boundary quadrilaterals, that are
triangulated by the couple of $2$-faces having $b$-sign equal to
$\pm$, respectively. The associated matrix dilogarithms are, by
definition, the QHFT tensors representing this bordism. They
represent the {\it elementary flip} on the quadrilateral
triangulations with two triangles. Any roughly flat/charged
$\Ii$-triangulated oriented pseudo-manifold $Z$ can be considered as
the result of a composition of several such tetrahedral elementary
bordisms. The fundamental transit invariance of the matrix
dilogarithms ensures {\it tautologically} that we have a well defined
so called {\it universal} QHFT for this very general bordism category.

\smallskip

The QHFTs constructed in the previous sections, which lead, in
particular, to the partition functions $\Hh_N(W,L_\Ff,\rho)$ or
$H_N(W,L,\rho)$ of Subsection \ref{PARTF}, naturally map into the
universal QHFT. Indeed, they are obtained via a specialization of its
setup (recall the further constraints of Def. \ref{flat/charge} for
flattenings and integral charges). The corresponding refined
flat/charged $\Ii$-equivalence classes {\it have a clear, intrinsic
  topological/geometric meaning}, described in Section 7 of
\cite{BB2}.

Another remarkable specialization of the universal QHFT environment is
given by the dilogarithmic invariants $H_N(M)$ of oriented non compact
complete hyperbolic $3$-manifolds $M$ of finite volume (for short:
{\it cusped manifolds}), defined in \cite {BB2}. These invariants are
constructed by using geodesic ideal triangulations with non negative
volume tetrahedra, which can be considered as triangulations of the
pseudo-manifold $\bar{M}$ obtained by taking the one point
compactification at each cuspidal end of $M$, where the vertices are
just the non-manifold points. Also in that situation there is a
natural notion of flattening and integral charge; the essential
difference is that for the integral charges we require only the global
constraint (a) in Definition \ref{flat/charge} (there is no link and
no boundary). By using a volume rigidity result for cusped
manifolds, we proved in \cite{BB2}, Th. 6.8, that the (normalized) flat/charged
$\Ii$-triangulations of a given cusped manifold $M$ all belong to the
same flat/charged $\Ii$-equivalence class in the universal QHFT
environment. Again, these classes have a clear, intrinsic
  topological/geometric meaning.

This is not evident for the equivalence class of an
arbitrary roughly flat/charged $\Ii$-triangulated pseudo-manifolds.
However, the above examples suggest the possibility that the universal QHFT
supports other {\it geometrically meaningful} specializations.

For instance, by using a suitable restricted (but non trivial) subset
of flat/charged $\Ii$-transits, we can formally construct a QHFT {\it
  free from the ``$=_N$'' phase ambiguity}. Clearly, the problem of
understanding the geometric meaning (if any) of this specialization of
the universal QHFT remains open. It is clearly related to the ultimate
understanding of the nature of the phase ambiguity itself.


 
\section*{Appendix}

\subsection*{A.1\ relationship between the QHFT marking of surfaces
and the Penner-Kashaev's coordinates}\label{COMPAREPAR}

As before, let $S$ be a closed compact surface of genus $g$ with a set
$V$ of $r$ marked points, and write $\hat{S}=S \setminus V$. Put a
finite area complete hyperbolic metric on $\hat{S}$. Denote by
$\mathcal{T}_g^r$ the marked Teichmuller space of $\hat{S}$. Recall
that $\mathcal{T}_g^r$ is homeomorphic to $\mr^{6g-6+2r}$, and to the
subspace of ${\rm Hom}(\pi_1(\hat{S}),PSL(2,\mr))/PSL(2,\mr)$ made of
{\it admissible} isomorphisms (i.e. mapping parabolic elements to
parabolics and induced by orientation preserving homotopy
equivalences), up to conjugation; this subspace lies in the manifold
part of ${\rm Hom}(\pi_1(\hat{S}),PSL(2,\mr))/PSL(2,\mr)$.

\smallskip

In \cite{Pe}, R.C. Penner constructed for any ideal triangulation $T'$
of $\hat{S}$ an $\mr_+^r$-principal bundle $\tilde{\Tt}_g^r
\rightarrow \Tt_g^r$ called the ``decorated'' Teichmuller space. The
total space $\tilde{\Tt}_g^r$ is parametrized by $-3\chi(\hat{S})$
coordinates associated to the edges of $T'$, called the {\it $\lambda$-lengths}. By fixing horocycles
about the punctures of $\hat{S}$ (more precisely, at the lifts to the
universal cover $\mathbb{D}^2$), the $\lambda$-lengths read as
$\sqrt{2\exp(\delta)}$, where $\delta$ is the algebraic distance between
the horocycles at the enpoints of the edges, counted positively if
they do not intersect. The radii of the horocycles form the semi-group
$\mr_+^r$, with the natural product action on $\tilde{\Tt}_g^r$. Recently, R.
Kashaev \cite{Kmod} generalized this construction to an $\mr_+^r$-bundle
$\tilde{\mathcal{M}} \rightarrow \mathcal{M}$, where $\mathcal{M}$ is
the space of conjugacy classes of irreducible
$PSL(2,\mr)$-representations of $\pi_1(\hat{S})$ with parabolic
holonomies around the punctures. Here we present a very simple way to
embed Kashaev's bundle $\tilde{\mathcal{M}}$ into the bundle of
cocycle $\Dd$-parameters of Definition \ref{cocpar}.

\paragraph{Kashaev's bundle $\tilde{\mathcal{M}}$.} Let us briefly
recall its ingredients. Given a conjugacy class of representations
$\rho\in \mathcal{M}$, one can find an ideal triangulation $T'$ of
$\hat{S}$ such that for any ideal arc of $T'$ the holonomies at the
endpoints lie in distinct parabolic subgroups of $PSL(2,\mr)$. By
removing an open small star at each vertex $v$ of $T'$, say star($v$),
we get a cellulation $C'$ of the surface $F$ introduced at the
beginning of Section \ref{Dpar}.  Each $2$-cell is a hexagon with $3$
`long' edges contained in one old edge of $T'$, and $3$ short edges
contained in the boundary link, link($v$) of star($v$).

\noindent Fix a
parabolic subgroup $P$ of $PSL(2,\mr)$, say upper-triangular, with
normalizer $B=N(P)$, the Borel subgroup of upper-triangular matrices.
We have the Bruhat decomposition $\textstyle PSL(2,\mr)=B \coprod
B\theta B$, where $\theta$ is the $(2\times 2)$-matrix with $1$ on the
bottom left, $-1$ on the top right, and $0$ elsewhere. Then, suitable
gauge transformations at the vertices of $C'$ allow us to construct
$1$-cocycles on $C'$ representing $\rho$ and such that the values
belong to the parabolic subgroup $P$ for all short edges, and to the
subset $\theta H$ for all long edges, where $H=B/P \cong \mr_+$. This
space of cocycles defines the fiber over $\rho$ in Kashaev's bundle
$\tilde{\mathcal{M}}$, and the residual gauge transformations turn out
to be isomorphic to a semi-group of $P$ isomorphic to $H^r \cong \mr_+^r$.

\noindent In fact, for each hexagon in $C'$
there exists a sign $\varepsilon=\pm 1$ such that the cocycle values
on the short edges are uniquely defined in terms of $\varepsilon$ and
the cocycle values on the long edges. So Kashaev's coordinates are
eventually given by:

\smallskip
 
-a sign $\varepsilon(t) \in \{-1,1\}$ for each triangle $t$ of $T'$
(i.e. each hexagon of $C'$);

-an $\mr_+$-valued function defined on the edges of $T'$ (i.e. the long edges of $C'$);

\smallskip

\noindent Kashaev described a partition of the fibration $\tilde{\mathcal{M}}
\rightarrow \mathcal{M}$ in terms of the above sign function
$\varepsilon$. The component where all the $\varepsilon$-values are positive or
negative is isomorphic to Penner's
decorated Teichmuller space $\tilde{\Tt}_g^r \rightarrow \Tt_g^r$.

\paragraph{Embedding $\tilde{\mathcal{M}}$ into $Z(T,b,1,y_1)_\beta$.}
By the two-dimensional version of \cite[Th. 3.4.9]{BP}, we can turn
$T'$ into an ideal triangulation of $F$ supporting a branching, via a
finite sequence of elementary flips. Fix a total ordering on the set
$V$ of marked points. By using the orientation of the short edges of
$C'$, an injective map $v \mapsto t_v$ as in Lemma \ref{injcorner}
selects one vertex $x_v$ of link($v$).  We stipulate that this is the
second endpoint of the short edge of $C'$ contained in $t_v$, and we
use $x_v$ as base point on link($v$).

\noindent Starting from
$x_v$, order the vertices of link($v$) in accordance with its
orientation. Then, in the reverse order, `slide' the vertices of
link($v$) one by one to $x_v$, except the biggest vertex $x_v''$, next
before $x_v$. This removes all the vertices on link($v$) but two.
After each sliding, the value of the cocycle on the new edge with one
endpoint at $x_v$ is forced. (The ordering on $V$ is needed because
each long edge shall eventually slide along two link($v$)'s). Do this
procedure for all the vertices of $T'$. We end up with a cellulation
of $F$ where each initial triangle $t_v
\subset T'$ appears as the union of one
triangle, one quadrilateral (with one short edge exactly), and one
bigon. We get also a well-determined, automatically defined,
$PSL(2,\mc)$-valued $1$-cocycle on this cellulation.

\noindent Remove a small monogon
inside each bigon and put a vertex $x_v'$ on its boundary, at the
'midpoint' of the segment $[x_v,x_v'']$. Triangulate the resulting
cellulation without adding new vertices. Clearly, we can extend the
branching of $T'$ during the above procedure. So we eventually find an
$e$-triangulation $(T,b)$ of $F$. At each step the cocycle values are
forced, except for the loop boundary edges, but also there the
holonomy about the loop is forced.

\noindent The constructions of
Section \ref{Dpar} imply that this cocycle can be turned into one of
the space $Z(T,b,1,y_1)_\beta$ via a suitable minimal sequence of gauge
transformations (recall that $1$ denotes the type of the
representations with only parabolic holonomies about the
punctures). Here we use the $PSL(2,\mr)$-version of the theory,
mentioned in Remark \ref{real}.
Also, the residual gauge transformations in Kashaev's bundle
$\tilde{\mathcal{M}}$ transit to the semi-group of $\Gg(T,b,1) \cong
{\rm Par}(2,\C)^r$ specified by the signs $c(g)=\pm 1$, for all
the boundary loop holonomies $g$.
So we eventually get an embedding $\tilde{\mathcal{M}}
\rightarrow Z(T,b,1,y_1)_\beta$.


\subsection*{A.2\ Formulas for the matrix dilogarithms}

Here we give the explicit formulas for the automorphisms
$\mathcal{R}_N(\Delta,b,w,f,c)$ associated in Section \ref{MATTENS} to flat/charged $\Ii$-tetrahedra, $N\geq 1$ being any
odd positive integer.
\smallskip

For $N=1$, we forget the integral charge $c$, so that
$\mathcal{R}_1$ is defined simply on flattened $\Ii$-tetrahedra. 
Namely, we have 
\begin{equation}\label{symRmatdil}
\mathcal{R}_1(\Delta,b,w,f) =
\exp\biggl(\frac{*_b}{i\pi}\biggl(-\frac{\pi^2}{6} - \frac{1}{2}
\int_0^{w_0} \biggl( \frac{{\rm l}_0(b,t,f)}{1-t} - \frac{{\rm
    l}_1(b,t,f)}{t} \biggr) \ dt \biggr) \biggr)
\end{equation}
where any ${\rm l}_j(b,t,f)$ is a log-branch
as defined in (\ref{logb1}). 
\medskip

For $N=2m+1>1$ and every complex number $x$ set
$x^{1/N} = \exp(\log(x)/N)$, where $\log$ is the standard branch of the
logarithm with arguments in $]-\pi,\pi]$
(by convention we put $0^{1/N} = 0$). 
Denote by $g$ the complex valued function, analytic
over the complex plane with cuts from the points $x =\zeta^k$ to
infinity ($k=1,\ldots,\ N-1$), defined by
$$g(x) := \prod_{j=1}^{N-1}(1 - x\zeta^{-j})^{j/N}$$
and set $h(x) := g(x)/g(1)$ (we have $\vert g(1) \vert =
N^{1/2}$). For any $u\in \C \setminus \{0,1\}$ and $ n \in \mathbb{N}$, put
$$\omega(u',v'\vert n) = \prod_{j=1}^n \frac{v'}{1-u'\zeta^j}$$ where
$v = 1-u$, and $u'$, $v'$ are arbitrary $N$th roots of $u$ and
$v$. These functions $\omega$ are periodic in their integer argument,
with period $N$. Finally, write $[x]=N^{-1}(1-x^N)/(1-x)$. Given a flat/charged $\Ii$-tetrahedron
$(\Delta,b,w,f,c)$, set
$$w_j'= \exp((1/N)(\log(w_j) + (f_j-*_bc_j)(N+1)\pi i)).$$
We define
\begin{equation}\label{symqmatdil}
\mathcal{R}_N(\Delta,b,w,f,c) =
\bigl((w_0')^{-c_1}(w_1')^{c_0}\bigr)^{\frac{N-1}{2}}\
(\Ll_N)^{*_b}(w_0',(w_1')^{-1})
\end{equation}
where 
$$\begin{array}{lll} \Ll_N(u',v')_{k,l}^{i,j}& = & h(u')\ \zeta^{kj+\frac{k^2}{2}}\
\omega(u',v'\vert i-k) \ \delta(i + j - l) \\
\bigl( \Ll_N(u',v')^{-1}\bigr)_{k,l}^{i,j}& = & \frac{[u']}{h(u')}\ \zeta^{-il-\frac{i^2}{2}}\
\frac{\delta(k + l - j)}{\omega(u'/\zeta,v'\vert k-i)}
\end{array}$$
and $\delta$ the $N$-periodic Kronecker symbol, i.e. $\delta(n) = 1$ if $n
\equiv 0$ mod($N$), and $\delta(n) = 0$ otherwise.



\end{document}